\documentclass[final,10pt]{elsarticle}
\usepackage{lineno,hyperref}
\usepackage{epsfig,amssymb,amsmath,version,amsthm}
\usepackage{amssymb,version,graphicx,fancybox,mathrsfs,multirow}
\usepackage[notcite,notref]{showkeys}
\usepackage{url,hyperref}
\usepackage{subfigure,epstopdf,bm}
\usepackage{color}
\usepackage[notcite,notref]{showkeys}
\usepackage{url,hyperref}
\usepackage{stmaryrd}
\usepackage{subfigure,epstopdf}
\usepackage{color}
\usepackage{cases}
\usepackage{accents}
\usepackage{appendix}
\usepackage{algorithm}
\usepackage{algorithmicx}
\usepackage{algpseudocode}
\usepackage{amsfonts}
\usepackage{booktabs} 
\usepackage{caption} 
\usepackage{makecell} 
\usepackage{graphicx, subcaption} 
\usepackage{graphicx} 
\usepackage{amsmath}   
\usepackage{amssymb} 

\modulolinenumbers[5]

\journal{Computer methods in applied mechanics and engineering}

\catcode`\@=11 \theoremstyle{plain}
\@addtoreset{equation}{section}

\@addtoreset{figure}{section}
\renewcommand\thefigure{\thesection.\@arabic\c@figure}
\@addtoreset{table}{section}
\renewcommand\thetable{\thesection.\@arabic\c@table}

\floatname{algorithm}{Algorithm}
\newtheorem{thm}{\bf Theorem}

\newtheorem{proposition}{Proposition}[section]

\newtheorem{cor}{\bf Corollary}

\newtheorem{lmm}{\bf Lemma}

\theoremstyle{remark}
\newtheorem{rem}{Remark}[section]

\def \ri {{\rm i}}
\newcommand{\bs}[1]{\boldsymbol{#1}}

\begin{document}

\begin{frontmatter}
\title{Fast boundary integral method for acoustic wave scattering in two-dimensional layered media}

\author[fn1]{Linfeng Xia}
\author[fn1]{Heng Yuan}
\author[fn1,fn2]{Bo Wang\corref{cor1}}
\ead{bowang@hunnu.edu.cn}
\author[fn3]{Wei Cai}

\cortext[cor1]{Corresponding author}

\address[fn1]{LCSM, Ministry of Education, School of Mathematics and Statistics, Hunan Normal University, Changsha, Hunan 410081, P. R. China.}
\address[fn2]{Xiangjiang Laboratory, Changsha, 410205, P. R. China.}
\address[fn3]{Department of Mathematics, Southern
Methodist University, Dallas, TX 75275, USA.}

\begin{abstract}
In this paper,  we present a  fast boundary integral method accelerated by the fast multipole method (FMM)  for acoustic wave scattering governed by the scalar Helmholtz equation in multi-layered two-dimensional media. Multiple scatterers are randomly distributed in the multi-layered medium with some scatterers possibly intersecting layer interfaces. The boundary integral formulation employs a layered-medium Green’s function to enforce transmission conditions across interfaces, thus avoiding unknowns on the interfaces and significantly reducing the size of the discretized problem  compared to approaches that use a free-space Green's function.  To demonstrate the FMM speedup,  a low-order Nystr\"{o}m method is used to discretize the boundary integral equation and then the resulting dense linear system is solved by GMRES iterative solver accelerated by an improved layered media FMM and a overlapping domain decomposition preconditioning. In the low frequency regime,  the proposed algorithm achieves an $\mathcal O(N)$ complexity. Numerical results  validate the accuracy, efficiency, and robustness of the method under complex settings and various incident angles. The proposed framework provides a scalable and efficient solver for acoustic wave scattering in layered media.
\end{abstract}


\end{frontmatter}

\section{Introduction}

Wave scattering in layered media is a fundamental problem in acoustics \cite{naqvi2025acoustic}, electromagnetics, and geophysics \cite{ha20253d}, with a wide range of applications from subsurface imaging to antenna design \cite{alibakhshikenari2021study} and optical waveguide analysis. The governing model is the Helmholtz equation posed in piecewise homogeneous domains separated by planar or curved interfaces, where continuity conditions across the interfaces must be enforced. Numerical simulations of such problems are notoriously challenging. Traditional domain-based discretization methods, such as the finite element method (FEM) or the finite difference method (FDM), require the introduction of artificial truncation boundaries to approximate the unbounded physical domain, typically combined with techniques such as perfectly matched layers (PML) \cite{chen2017pml,jiang2019adaptive, dastour2019fourth,jiang2022pml} or absorbing boundary conditions (ABCs) \cite{Absorb1977,higdon1992absorbing,aksun2002comparative,wu2021high}. While these approaches are highly versatile, they generally lead to very large linear systems, particularly in high-frequency regimes, rendering large-scale computations prohibitively expensive. Moreover, the construction of the ABCs and the stability of the PML technique for multi-layered media are still challenging. 

Boundary integral equation (BIE) methods \cite{kress1991boundary,bendali2008boundary,hsiao2021boundary,zaman2000comprehensive} provide an attractive alternative to address these challenges. By using Green’s functions, BIE reformulates the governing partial differential equations into integral equations posed solely on the boundary of the scatterer. This approach offers  significant advantages of dimensional reduction, which greatly reduces the number of unknowns. Furthermore, BIE methods inherently enforce the radiation condition at infinity, eliminating the need for artificial domain truncation and thus avoiding the need of PML-type or ABCs techniques. These properties make BIE methods especially suitable for exterior scattering problems.

Despite these intrinsic advantages, the application of BIEs to layered media with complex and multiscale structures encounters two principal computational bottlenecks. First, the Nystr\"{o}m discretization results in dense and typically non-symmetric linear system. For a problem with $N$ unknowns, classical direct solvers (e.g., Gaussian elimination) entail both storage and computational complexities of order $O(N^3)$, which is prohibitively expensive for problems with large numbers of unknowns. Second, the evaluation of layered Green’s functions is itself highly challenging, as they are represented by Sommerfeld integrals \cite{sommerfeld1909ausbreitung,mosig2021sommerfeld} with oscillatory and slowly decaying kernels, often requiring the numerical computation of semi-infinite integrals significantly more expensive than the evaluation of free-space Green’s functions.

To address the first bottleneck associated with dense linear system, a broad spectrum of fast algorithms has been developed over the past three decades. For instance, Hackbusch and collaborators introduced the hierarchical matrix ($\mathcal{H}$-matrix) framework \cite{wolfgang2015hierarchical}, which leverages block-wise low-rank approximations to construct data-sparse representations, enabling fast matrix–vector multiplications, approximate LU factorization, and direct solvers. Bebendorf and co-workers further developed adaptive cross approximation (ACA) and block low-rank methods \cite{bebendorf2000approximation}, which construct low-rank representations adaptively from a small number of row and column samples, significantly reducing assembly and storage costs. Additional progress includes Fourier integral operator (FIO)-oriented methods suitable for high-frequency regimes \cite{candes2009fast}, as well as recursive skeletonization and other fast direct solvers designed for problems with multiple right-hand sides or highly ill-conditioned systems \cite{martinsson2005fast, ho2012fast}. Among these, the fast multipole method (FMM) introduced by Greengard and Rokhlin \cite{greengard1987fast} stands out for its remarkable efficiency: by hierarchically compressing far-field interactions via multipole and local expansions, it reduces the cost of matrix–vector multiplications from $O(N^2)$ to nearly linear complexity ($O(N)$ or $O(N \log N)$). Subsequent extensions have adapted FMM to BIEs \cite{ying2009fast} and layered media FMMs \cite{wang2019fast}, thereby retaining scalability while effectively accounting for interface reflections and transmissions.

Applying these fast algorithms to layered media scattering problems introduces the second major bottleneck: the expensive evaluation of layered Green’s functions. Early techniques include the discrete complex image method (DCIM) and contour deformation approaches, but both face substantial limitations for multilayered configurations. To address these challenges, Cai and Yu developed a window-based acceleration strategy \cite{cai2000fast} to mitigate the slow convergence of the Sommerfeld integrals caused by surface-wave poles, the slow decay of spectral Green’s functions, and oscillations in Hankel-transform kernels. Michalski and Mosig later provided a comprehensive review of efficient evaluation strategies for the Sommerfeld integrals \cite{michalski2016efficient}, focusing on handling the oscillatory and singular behaviors of layered kernels. More recent advances have introduced even more efficient and robust methods, including the windowed Green’s function (WGF) method \cite{Bruno2016_WGF}, the PML-BIE method \cite{bao2024highly,wang2024fft,lu2023highly}, hybrid integral formulations \cite{lai2018new}, and numerical quadratures of high-accuracy such as the Tanh–Sinh quadrature for Sommerfeld-type integrals \cite{vanherck2020tanh}. 

In this work, we investigate the two-dimensional acoustic wave scattering problem in multilayered media and propose a fast boundary integral method accelerated by FMM for layered Green's function of the Helmholtz equation. First, We derive the layered Green’s function that satisfies the transmission conditions across all interfaces, which allows the resulting boundary integral formulation to substantially reduce the number of unknowns. Secondly,  building on the FMM developed for the Helmholtz equation in layered media \cite{zhang2020exponential,wang2019fast}, we improve the two-dimensional setting with an enhancement of equivalent polarization coordinate and effective location techniques proposed in \cite{Yuan2025Fast}. Specifically, we give an more accurate description of the vertical transmission distance of the reaction-field components and employ polarization coordinates for sources as well as effective locations for targets, thereby improving both the accuracy and efficiency of the FMM. With FMM acceleration, the solution of linear systems arising from low-order Nystr\"{o}m discretizations achieves nearly linear computational complexity, enabling efficient simulation of large-scale problems. Furthermore, inspired by the preconditioner based on the hierarchical tree structure of the FMM (cf. \cite{nabors2002fast}), we develop an effective preconditioning strategy that mitigates the deterioration of the system’s condition number as the problem size increases, reducing both the iteration count and the computational time by roughly one half. Extensive numerical experiments demonstrate that the proposed method is accurate, robust, and highly efficient, achieving nearly linear computational complexity while maintaining first-order numerical accuracy.

The remainder of the paper is organized as follows.
In Section \ref{sect2}, we formulate the two-dimensional Helmholtz scattering problem in a multi-layered medium with a plane incident wave from the top layer. By splitting out the background scattering field, an interface problem for the purely outgoing wave is formulated. Then, a boundary integral method with Nystr\"{o}m dicretization is presented 
in Section \ref{sect3}. Here, we use the layered Green’s function to naturally impose transmission conditions without introducing redundant unknowns on the interface. 
In Section \ref{sect4}, an improved fast multipole algorithm  is proposed to accelerate the matrix-vector product in the GMRES iterative algorithm and an overlapping domain decomposition preconditioner is presented based on the tree structure of the FMM to reduce the number of iterations. 
In Section \ref{sect5}, numerical examples with multiple scatterers randomly distributed in multi-layered media are provided to validate the accuracy and efficiency of the proposed algorithm. 

\section{Acoustic wave scattering problem in two-dimensional layered media}\label{sect2}
\subsection{Model problem}
\begin{figure}[htbp]
	\centering
	\includegraphics[scale=0.35]{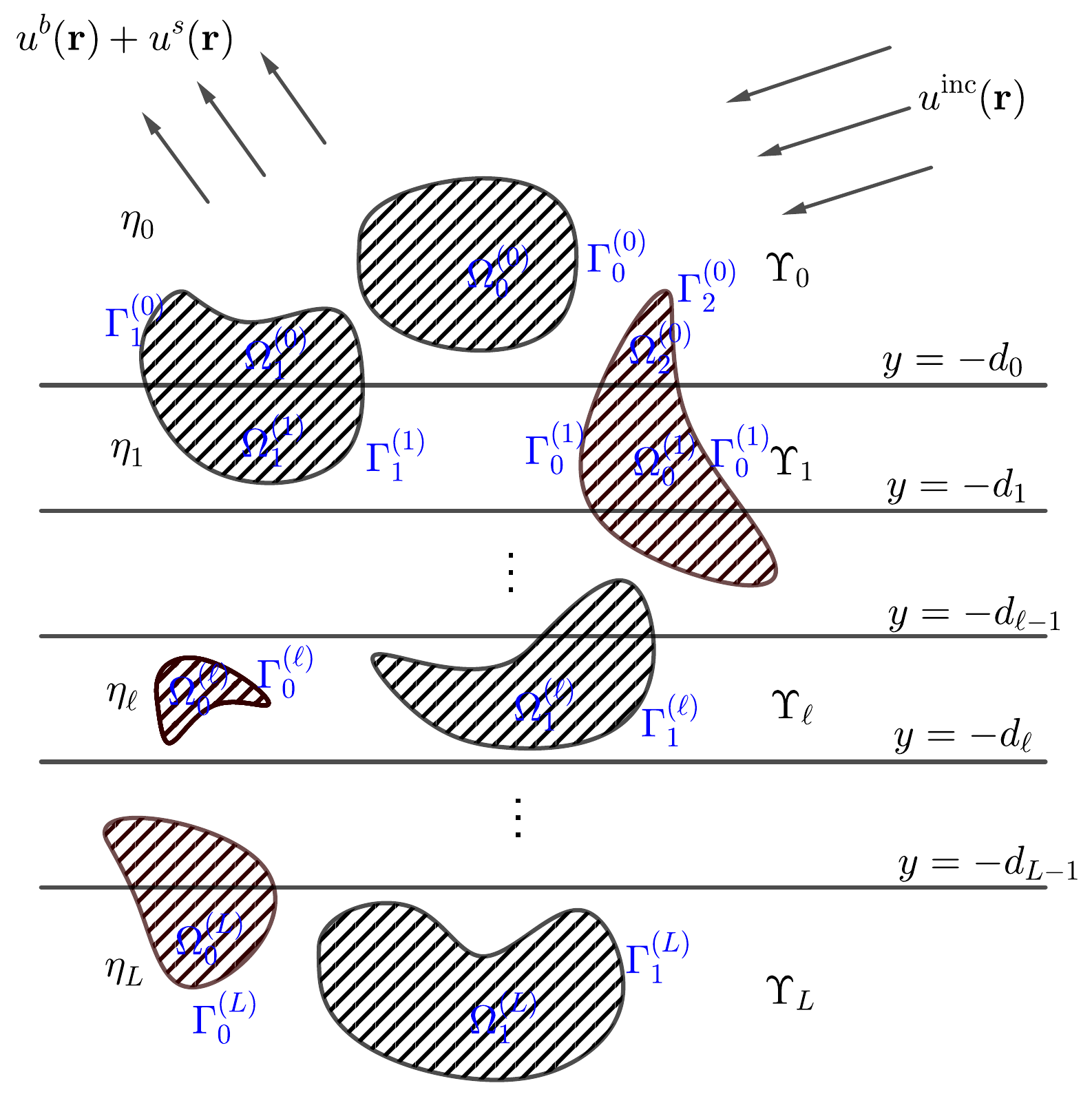}
	\caption{Configurations of the scattering problem in layered media.}
	\label{scatterersconfigure}
\end{figure}

Consider the layered media of $L+1$ layers with interfaces $y=-d_{\ell}$, $\ell=0, 1, \cdots, L-1$ and reflection index $\{\eta_{\ell}\}_{\ell=0}^{L}$. Let $\{D_{j}\}_{j=1}^M$ be $M$ impenetrable scatterers randomly distributed in the layered media and $\{\partial D_{j}\}_{j=1}^M$ define the boundaries of scatterers $\{D_{j}\}_{j=1}^M$. Each of the scatterers can be located in one layer or cross multiple layers, see Fig. \ref{scatterersconfigure} . Denote by 
\begin{equation}
	\begin{split}
		\Upsilon_0=&\{(x, y)| y>-d_0\},\quad \Upsilon_{L}=\{(x, y)| y<-d_{L-1}\},\\ 
		\Upsilon_{\ell}=&\{(x, y)| -d_{\ell}<y<-d_{\ell-1}\},\quad \ell=1, 2, \cdots, L-1,
	\end{split}
\end{equation}
the infinite strips of the layers, 
\begin{equation}
	\Omega_{\ell}=\bigcup_{i=1}^M\Omega_{i}^{(\ell)}=  \bigcup_{i=1}^M(\Upsilon_{\ell}\cap D_i),\quad\ell=0, 1, \cdots, L,
\end{equation}
the disconnected domains in the $\ell$-th layer occupied by the multiple scatterers, and 
\begin{equation}
    \Gamma_{\ell} = \bigcup_{i=1}^M\Gamma_{i}^{(\ell)}=  \bigcup_{i=1}^M(\Upsilon_{\ell}\cap \partial D_i),\quad\ell=0, 1, \cdots, L,
\end{equation}
the boundary of the disconnected domain.

The mathematical model of a time-harmonic acoustic scattering problem under this general setting  is the Helmholtz equation of the form 
\begin{equation}\label{TLhelmholtztotallayer1}
		(\bs{\Delta}+k^2_{\ell})u(\bs r)=0, \quad \bs r\in\Upsilon_{\ell}\backslash\Omega_{\ell},\quad \ell=0, 1, \cdots, L,
\end{equation}
in each layer. Here, $u=u^{\rm sc}(\bs r)+u^{\rm inc}(\bs r)$ is the total field, $k_{\ell}$ is the wave number in the $\ell$-th layer. The incident wave $u^{\rm inc}(\bs r)$ is typically a plane wave given in the $0$-th layer. Its expression is given by
\begin{equation}\label{incidentwave}
	u^{\rm inc}(\bs r)=\begin{cases}
	\displaystyle	e^{\ri (k_{0,x}x-k_{0,y}y)}, &  y>-d_0,\\
	\displaystyle 0, & {\rm otherwise},
	\end{cases}
\end{equation}
where $(k_{0,x}, k_{0,y})$ is the propagating vector such that $k_{0,y}>0$, and $(k_{0,x})^2+(k_{0,y})^2=k_0^2$. 
Due to the presence of the layered media, transmission conditions 
\begin{equation}
	\llbracket u\rrbracket=0,\; \Big\llbracket \eta \frac{\partial u}{\partial\bs n}\Big\rrbracket=0, \quad y=-d_{\ell}, \;\; \ell =0,1,\cdots,L-1,
\end{equation}
are imposed on the interfaces. On the boundaries of the scatterers, the Dirichlet, Neumann or Robin boundary conditions could be imposed according to different  materials of the scatterers. For simplicity of presentation, only sound soft (Dirichlet) boundary condition,
i.e, 
\begin{equation}
	u(\bs r)=0,\quad \bs r\in\Omega_{\ell},\quad \ell=0, 1, \cdots, L
\end{equation}
will be considered in this paper.

\subsection{Background scattering field}
The reflection of plane waves by a layered structure can be calculated analytically. In order to formulate BIEs for the scattering problem, we separate the background scattering field from the scattering field $u^{sc}(\bm{r})$. The background scattering field denoted by $u^b(\bs r)$ is the reflected fields within the layered media in the absence of scatterers.  It satisfies
\begin{equation}
\begin{split}
\label{TLhelmholtztotallayer2}
(\bs{\Delta}+k^2_{\ell})u^b(\bs r)=0, &\quad \bs r\in\Upsilon_{\ell},\quad \ell=0, 1, \cdots, L,\\
\llbracket u^b+u^{\rm inc}\rrbracket=0, \quad
\Big\llbracket \eta\Big( \frac{\partial u^b}{\partial\bs n}+ \frac{\partial u^{\rm inc}}{\partial\bs n}\Big)\Big\rrbracket=0, &\quad y=-d_{0},\\
\llbracket u^b\rrbracket=0, \quad
\Big\llbracket \eta\frac{\partial u^b}{\partial\bs n}\Big\rrbracket=0, &\quad y=-d_{\ell}, \;\; \ell =1,\cdots,L-1,
\end{split}
\end{equation}
it is well-known that $u^b(\bs r)$ has an ansatz expression
\begin{equation}
	u^b(\bs r)=\begin{cases}
		A_{0}e^{\ri(k_{0,x}x+k_{0,y}y)},\quad     \bs r\in\Upsilon_{0},\\
        A_{\ell}e^{\ri(k_{\ell,x}x+k_{\ell,y}y)}+B_{\ell} e^{\ri(k_{\ell,x}x-k_{\ell,y}y)},\quad \bs r\in\Upsilon_{\ell},\;\;\ell\neq 0, L,\\
		B_{L} e^{\ri(k_{L,x}x-k_{L,y}y)},\quad \bs r\in\Upsilon_{L}.
	\end{cases}
\end{equation}
Applying the interface conditions gives
\begin{equation}
	\begin{split}
		k_{L,x}=&k_{L-1,x}=\cdots=k_{1,x}=k_{0,x},\quad k_{\ell,y}= \sqrt{{k^2_{\ell}}-{k^2_{\ell,x}}},\\ 
		A_{0}=&\widetilde R_{01}e^{2{\rm i} k_{0,y}d_0},\quad B_{L}=\widetilde T_{0L},\\
		A_{\ell}=&\widetilde R_{\ell\ell+1}\widetilde T_{0\ell}e^{2{\rm i} k_{\ell,y}d_{\ell}},\quad B_{\ell}=\widetilde T_{0\ell},\quad \ell=1, 2, \cdots, L-1,
	\end{split}
\end{equation}
where the general reflection coefficients $\{\widetilde R_{\ell\ell+1}(\lambda)\}_{\ell=0}^{L-1}$ and general transmission coefficients $\{\widetilde T_{0\ell}(\lambda)\}_{\ell=1}^{L}$ are defined in \eqref{Rellrecursiondown} and \eqref{generaltransmission}, respectively. 

Then, the scattering field due to the presence of the scatterers is defined as
$$
u^s(\bs r)=u(\bs r)-u^b(\bs r)-u^{\rm inc}(\bs r),
$$
which is the solution of the following scattering problem
\begin{equation}\label{extdirichleteq1}
    \begin{aligned}
        (\bs{\Delta}+k^2_{\ell})u^s(\bs r)=0, &\quad \bs r\in\Upsilon_{\ell}\setminus\Omega_{\ell},\quad \ell=0, 1, \cdots, L,\\
        u^s(\bs r)=g(\bs r),&\quad \bs r\in\Gamma_{\ell},\quad \ell=0, 1, \cdots, L,\\
        \llbracket u^s\rrbracket=0,\quad \Big\llbracket \eta \frac{\partial u^s}{\partial\bs n}\Big\rrbracket=0,&\quad y=-d_{\ell},\quad \ell=0, 1, \cdots, L-1,\\
        u^s(\bs r)\;{\rm is \; outing\; as} & \quad |\bm r|\rightarrow\infty.
    \end{aligned}
\end{equation}
Here, the boundary data is given by
\begin{equation}\label{homogeneousscattering}
	g(\bs r)=\begin{cases}
		\displaystyle	-u^{\rm inc}(\bs r)-u^b(\bs r),&\quad \bs r\in\Gamma_0,\\
		\displaystyle	-u^b(\bs r),&\quad \bs r\in\bigcup_{\ell=1}^{L}\Gamma_{\ell},\\
	\end{cases}
\end{equation}
which implies that the actual incident wave is the superposition of the original incident wave and the background scattering wave. 
In the rest of this paper, we will employ the boundary integral technique to develop a fast algorithm to solve the wave scattering problem \eqref{extdirichleteq1}.

\section{Boundary integral equation and Nystr\"{o}m discretization}\label{sect3}
In this section, we employ the layered Green's function \eqref{layeredGreensfun} to derive BIEs for the model problem \eqref{extdirichleteq1} and then develop a fast numerical algorithm based on Nystr\"{o}m discretization and FMM accelerated GMRES solver. Due to the use of layered Green's function, The transmission interface conditions are automatically satisfied by the integral representation of the solution and the integral equation is formulated solely on the boundaries of the scatterers.

\subsection{Green's function of the Helmholtz equation in layered media} 
Suppose we have a point
source at $\boldsymbol{r}^{\prime}=(x^{\prime},y^{\prime})$ in the
$\ell^{\prime}$th layer ($-d_{\ell^{\prime}}<y^{\prime}<-d_{\ell^{\prime}-1}$) of the layered medium depicted in Fig. \ref{scatterersconfigure}, the layered media Green's function $u_{\ell\ell'}(\bs r, \bs r')$ for Helmholtz equation satisfies
\begin{equation}\label{MutiLayerHelmEquation}
	(\bs{\Delta}+k^2_{\ell})G_{\ell\ell'}(\boldsymbol{r},\boldsymbol{r}^{\prime
	})=-\delta(\boldsymbol{r},\boldsymbol{r}^{\prime}),
\end{equation}
at field point $\boldsymbol{r}=(x,y)$ in the $\ell$th layer ($-d_{\ell}<y<-d_{\ell}-1$) where $\delta(\boldsymbol{r},\boldsymbol{r}^{\prime})$ is the
Dirac delta function.

For any $\ell=0 \to L-1$, the reflection and transmission coefficients due to the interface $y=-d_{\ell}$ are given by
\begin{equation}
	\begin{split}
		R_{\ell\ell+1}(\lambda)=\frac{\eta_{\ell}k_{\ell, y} - \eta_{\ell+1}k_{\ell+1, y}}{\eta_{\ell}k_{\ell, y}+\eta_{\ell+1}k_{\ell+1, y}},\quad T_{\ell\ell+1}(\lambda)=\frac{2\eta_{\ell}k_{\ell, y}}{\eta_{\ell}k_{\ell, y} + \eta_{\ell+1}k_{\ell+1, y}},
	\end{split}
\end{equation}
if a plane wave $e^{\ri(\lambda x+k_{\ell,y}y)}$ is incident from above. In contrast, if a plane wave $e^{\ri (\lambda x+k_{\ell+1, y}y)}$ is incient from the $(\ell+1)$-th layer from below,  the reflection and transmission coefficients are 
\begin{equation}
	R_{\ell+1,\ell}=-R_{\ell\ell+1},\quad T_{\ell+1,\ell}=\frac{\eta_{\ell+1}k_{\ell+1, y}}{\eta_{\ell}k_{\ell, y}}T_{\ell\ell+1}.
\end{equation}
. 

In order to investigate the multi-layered problem, the general reflection and transmission coefficients were introduced (cf. \cite{chew1999waves}). Two groups of general reflection coefficients $\{\widetilde R_{\ell+1,\ell}(\lambda)\}_{\ell=0}^{L-1}$, $\{\widetilde R_{\ell\ell+1}(\lambda)\}_{\ell=0}^{L-1}$ can be calculated via two recursions as follows
\begin{itemize}
	\item Set $\widetilde R_{0,-1}(\lambda)=0$, then
	\begin{equation}\label{Rellrecursionup}
		\widetilde R_{\ell+1,\ell}(\lambda)
		=\frac{R_{\ell+1,\ell}+\widetilde R_{\ell\ell-1}(\lambda)e^{-2 k_{\ell,y}(d_{\ell}-d_{\ell-1})}}{1+R_{\ell+1,\ell}\widetilde R_{\ell\ell-1}(\lambda)e^{-2 k_{\ell,y}(d_{\ell}-d_{\ell-1})}},\quad \ell=0, 1, \cdots, L-1.
	\end{equation}
	\item Set $\widetilde R_{L-1,L}(\lambda)=0$, then
	\begin{equation}\label{Rellrecursiondown}
		\widetilde R_{\ell\ell+1}(\lambda)
		=\frac{R_{\ell\ell+1}+\widetilde R_{\ell+1,\ell+2}(\lambda)e^{-2k_{\ell+1,y}(d_{\ell+1}-d_{\ell})}}{1+R_{\ell\ell+1}\widetilde R_{\ell+1,\ell+2}(\lambda)e^{-2 k_{\ell+1,y}(d_{\ell+1}-d_{\ell})}},\quad \ell=L-2, L-3, \cdots, 0.
	\end{equation}
\end{itemize}
The general transmission coefficient from $\ell'$-th layer to $\ell$-th layer has recursion
\begin{equation}\label{generaltransmission}
	\widetilde T_{\ell'\ell}(\lambda)=\begin{cases}
		\displaystyle	\frac{T_{\ell+1,\ell}(\lambda)e^{\ri( k_{\ell,y}-k_{\ell+1,y})d_{\ell}}\widetilde T_{\ell'\ell+1}(\lambda)}{1+R_{\ell+1,\ell}(\lambda)\widetilde R_{\ell,\ell-1}(\lambda)e^{-2 k_{\ell,y}(d_{\ell}-d_{\ell-1})}},& \ell=\ell'-1,\ell'-2, \cdots, 0,\\
		\displaystyle\frac{T_{\ell-1,\ell}(\lambda)e^{\ri (k_{\ell-1,y}-k_{\ell,y})d_{\ell-1}}\widetilde T_{\ell'\ell-1}(\lambda)}{1+R_{\ell-1,\ell}(\lambda)\widetilde R_{\ell,\ell+1}(\lambda)e^{-2k_{\ell,y}(d_{\ell}-d_{\ell-1})}},& \ell=\ell'+1, \ell'+2, \cdots, L,
	\end{cases}
\end{equation}
where the initial value is given by $\widetilde T_{\ell'\ell'}(\lambda)\equiv 1$.

By using Fourier transforms along $x-$directions, the problem (\ref{MutiLayerHelmEquation}) can be solved analytically for each layer in $y$ by imposing
transmission conditions at the interface between $\ell$th and $(\ell-1)$th
layer ($y=-d_{\ell-1})$, \textit{i.e.},
\begin{equation}
	G_{\ell-1,\ell'}(\bm{r},\bm{r}')=G_{\ell\ell'}(\bm{r},\bm{r}'),\quad \eta_{\ell-1}\frac{\partial  G_{\ell-1,\ell'}(\bm{r},\bm{r}')}{\partial y}=\eta_{\ell}\frac{\partial G_{\ell\ell'}(\bm{r},\bm{r}')}{\partial y},
\end{equation}
as well as the outgoing conditions in the top and bottom-most layers as
$y\rightarrow\pm\infty$.
Here, we just present the expression of the layered Green's function. Detailed derivation can be found in (cf. \cite{chew1999waves,wang2020taylor}). The expression of the Green's function in the physical domain takes the form
\begin{equation}\label{layeredGreensfun}
	G_{\ell\ell^{\prime}}(\boldsymbol{r},\boldsymbol{r}^{\prime})=\delta_{\ell\ell'}G^f_{\ell}(\bs r, \bs r')+G^{\rm r}_{\ell\ell'}(\boldsymbol{r},\boldsymbol{r}^{\prime}),
\end{equation}
where $\delta_{\ell\ell'}$ is the Kronecker symbol, 
\begin{equation}\label{freegreensfun}
G^f_{\ell}(\bs r, \bs r')=\frac{\ri}{4}H_0^{(1)}(k_{\ell}|\bs r-\bs r^{\prime}|)
\end{equation}
is the free space Green's function with wave number $k_{\ell}$, $G^{\rm r}_{\ell\ell'}(\bs r, \bs r')$ is the reaction component. In general, $G^{\rm r}_{\ell\ell'}(\bs r, \bs r')$ does not have closed form. It is usually given by the inverse Fourier transform of the Green's function in the spectral domain. According to the dependence on the $y$ and $y'$ coordinates, the reaction component $G^{\rm r}_{\ell\ell'}(\bs r, \bs r')$ generally consists of four components, i.e., 
\begin{equation}\label{reactioncomponent}
	G^{\rm r}_{\ell\ell'}(\boldsymbol{r},\boldsymbol{r}^{\prime})=G_{\ell\ell^{\prime}}^{\uparrow\downarrow}(\boldsymbol{r},\boldsymbol{r}^{\prime})+G_{\ell\ell^{\prime}}^{\uparrow\uparrow}(\boldsymbol{r},\boldsymbol{r}^{\prime})+G_{\ell\ell^{\prime}}^{\downarrow\downarrow}(\boldsymbol{r},\boldsymbol{r}^{\prime})+G_{\ell\ell^{\prime}}^{\downarrow\uparrow}(\boldsymbol{r},\boldsymbol{r}^{\prime}),
\end{equation}
where 
\begin{equation}\label{fourreactions}
	\begin{split}
		G_{\ell\ell^{\prime}}^{\uparrow\downarrow}(\bs r, \bs r')=&\begin{cases}
			\displaystyle \frac{1}{2\pi }\int_{-\infty}^{\infty}\frac{ e^{\ri \lambda(x-x')+\ri(k_{\ell,y}y-k_{\ell',y}\tau_{\ell'}(y'))}}{2k_{\ell',y}}\sigma^{\uparrow\downarrow}_{\ell\ell'}(\lambda) d\lambda,\quad 0\leq\ell\leq \ell',\\
			\displaystyle \frac{1}{2\pi }\int_{-\infty}^{\infty}\frac{ e^{\ri \lambda(x-x')+\ri (k_{\ell',y}y'-k_{\ell,y}\tau_{\ell}(y))}}{2k_{\ell',y}}\sigma^{\uparrow\downarrow}_{\ell\ell'}(\lambda) d\lambda,\quad \ell'<\ell\leq L,
		\end{cases}\\
		G_{\ell\ell^{\prime}}^{\uparrow\uparrow}(\bs r, \bs r')=&\begin{cases}
			\displaystyle \frac{1}{2\pi }\int_{-\infty}^{\infty}\frac{ e^{\ri \lambda(x-x')+ \ri(k_{\ell,y}y-k_{\ell',y}y')}}{2k_{\ell',y}}\sigma^{\uparrow\uparrow}_{\ell \ell'}(\lambda)  d\lambda,\quad 0\leq\ell<\ell',\\
			\displaystyle \frac{1}{2\pi }\int_{-\infty}^{\infty}\frac{ e^{\ri \lambda(x-x')+\ri (k_{\ell',y}\tau_{\ell'-1}(y')-k_{\ell,y}\tau_{\ell}(y))}}{2k_{\ell',y}}\sigma^{\uparrow\uparrow}_{\ell\ell'}(\lambda) d\lambda,\quad \ell'\leq\ell\leq L,
		\end{cases}\\
		G_{\ell\ell^{\prime}}^{\downarrow\downarrow}(\bs r, \bs r')=&\begin{cases}
			\displaystyle \frac{1}{2\pi }\int_{-\infty}^{\infty}\frac{ e^{\ri \lambda(x-x')+\ri (k_{\ell,y}\tau_{\ell-1}(y)-k_{\ell',y}\tau_{\ell'}(y'))}}{2k_{\ell',y}}\sigma^{\downarrow\downarrow}_{\ell\ell'}(\lambda) d\lambda,\quad 0\leq\ell\leq \ell',\\	
			\displaystyle \frac{1}{2\pi }\int_{-\infty}^{\infty}\frac{ e^{\ri \lambda(x-x')+\ri (k_{\ell',y}y'-k_{\ell,y}y)}}{2k_{\ell',y}}\sigma^{\downarrow\downarrow}_{\ell\ell'}(\lambda) d\lambda,\quad \ell'<\ell\leq L,
		\end{cases}\\
		G_{\ell\ell^{\prime}}^{\downarrow\uparrow}(\bs r, \bs r')=&\begin{cases}
			\displaystyle	\frac{1}{2\pi }\int_{-\infty}^{\infty}\frac{ e^{\ri \lambda(x-x')+\ri (k_{\ell,y}\tau_{\ell-1}(y)-k_{\ell',y}y')}}{2k_{\ell',y}}\sigma^{\downarrow\uparrow}_{\ell \ell'}(\lambda) d\lambda,\quad 0\leq \ell<\ell',\\
			\displaystyle	\frac{1}{2\pi }\int_{-\infty}^{\infty}\frac{ e^{\ri \lambda(x-x')+\ri (k_{\ell',y}\tau_{\ell'-1}(y')-k_{\ell,y}y)}}{2k_{\ell',y}}\sigma^{\downarrow\uparrow}_{\ell\ell'}(\lambda) d\lambda,\quad \ell'\leq\ell\leq L.
		\end{cases}
	\end{split}
\end{equation}
In the above integral representations, 
\begin{equation}\label{reflections}
	\tau_{m}(a)=-2d_m-a,\quad m=0, 1, \cdots, L-1
\end{equation}
are the reflection of any given $a$ according to the interface $y=-d_m$, and the densities are given by
\begin{itemize}
		\item {\bf Case I:} $\ell=\ell'$,
		\begin{equation}\label{densityinsourcelayer}
			\begin{split}
				&\sigma_{\ell'\ell'}^{\uparrow\downarrow}(\lambda)=\frac{\widetilde R_{\ell'\ell'+1}}{1-\widetilde R_{\ell\ell+1}\widetilde R_{\ell\ell-1}e^{-2 k_{\ell',y}(d_{\ell'}-d_{\ell'-1})}},\\ &\sigma_{\ell'\ell'}^{\downarrow\uparrow}(\lambda)=\frac{\widetilde R_{\ell'\ell'-1}}{1-\widetilde R_{\ell\ell+1}\widetilde R_{\ell\ell-1}e^{-2 k_{\ell',y}(d_{\ell'}-d_{\ell'-1})}},\\
				&\sigma_{\ell'\ell'}^{\uparrow\uparrow}(\lambda)=\sigma_{\ell'\ell'}^{\downarrow\downarrow}(\lambda)=\widetilde R_{\ell'\ell'-1}\sigma_{\ell\ell'}^{\uparrow\downarrow}(\lambda)=\widetilde R_{\ell'\ell'+1}\sigma_{\ell\ell'}^{\downarrow\uparrow}(\lambda).
			\end{split}
		\end{equation}
		\item {\bf Case II:} $0\leq \ell<\ell'$,
		\begin{equation}\label{disity1}
			\begin{split}
				& \sigma^{\uparrow\downarrow}_{\ell \ell'}(\lambda)=\widetilde T_{\ell'\ell}(\lambda)\sigma^{\uparrow\downarrow}_{\ell'\ell'}(\lambda), \quad \sigma^{\uparrow\uparrow}_{\ell \ell'}(\lambda)=\widetilde T_{\ell'\ell}(\lambda)\Big[1+\sigma_{\ell'\ell'}^{\uparrow\uparrow}(\lambda)e^{2{\rm{i}}k_{\ell',y}(d_{\ell'}-d_{\ell'-1})}\Big], \\
				& \sigma^{\downarrow\downarrow}_{\ell \ell'}(\lambda) = \widetilde R_{\ell\ell-1}(\lambda)\sigma^{\uparrow\downarrow}_{\ell \ell'}(\lambda), \quad \sigma^{\downarrow\uparrow}_{\ell \ell'}(\lambda)=\widetilde R_{\ell\ell-1}(\lambda)\sigma^{\uparrow\uparrow}_{\ell \ell'}(\lambda).
			\end{split}
		\end{equation}
		\item {\bf Case III:} $\ell'< \ell\leq L$
		\begin{equation}\label{disity2}
			\begin{split}
				& \sigma^{\downarrow\uparrow}_{\ell\ell'}(\lambda)=\widetilde T_{\ell'\ell}(\lambda)\sigma^{\downarrow\uparrow}_{\ell'\ell'}(\lambda), \quad \sigma^{\downarrow\downarrow}_{\ell\ell'}(\lambda)=\widetilde T_{\ell'\ell}(\lambda)\Big[1+\sigma_{\ell'\ell'}^{\downarrow\downarrow}(\lambda)e^{2{\rm{i}}k_{\ell',y}(d_{\ell'}-d_{\ell'-1})}\Big], \\
				& \sigma^{\uparrow\uparrow}_{\ell\ell'}(\lambda)=\widetilde R_{\ell\ell+1}(\lambda)\sigma^{\downarrow\uparrow}_{\ell\ell'}(\lambda), \quad \sigma^{\uparrow\downarrow}_{\ell\ell'}(\lambda)=\widetilde R_{\ell\ell+1}(\lambda)\sigma^{\downarrow\downarrow}_{\ell\ell'}(\lambda).
			\end{split}
		\end{equation}
\end{itemize}

As in our previous work \cite{Yuan2025Fast}, we introduce the \emph{equivalent polarization coordinates}
\begin{equation}\label{equivalent_imag_for_source}
\breve\tau_{\ell\ell'}^{\downarrow}(\bs r')=\begin{cases}
			\tau_{\ell'}(\bs r'),& \ell\le\ell',\\
			\bs r',& \ell>\ell',
		\end{cases}\quad
		\breve\tau_{\ell\ell'}^{\uparrow}(\bs r')=\begin{cases}
			\bs r',&\ell<\ell',\\
			\tau_{\ell'-1}(\bs r'),&\ell\ge\ell',
		\end{cases}
\end{equation}
for source $\bs r'$ and \emph{effective locations}
\begin{equation}\label{equivalent_imag_for_target}
\begin{split}
&\hat\tau_{\ell\ell'}^{\uparrow\uparrow}(\bs r)=\begin{cases}
\bs r,& \ell<\ell',\\
\tau_{\ell}(\bs r),& \ell\ge\ell',
\end{cases}\quad\hat\tau_{\ell\ell'}^{\uparrow\downarrow}(\bs r)=\begin{cases}
\bs r,& \ell\le\ell',\\
\tau_{\ell}(\bs r),& \ell> \ell',
\end{cases}\\
&\hat\tau_{\ell\ell'}^{\downarrow\uparrow}(\bs r)=\begin{cases}
\tau_{\ell-1}(\bs r),& \ell<\ell',\\
\bs r,& \ell\ge\ell',
\end{cases}\quad \hat\tau_{\ell\ell'}^{\downarrow\downarrow}(\bs r)=\begin{cases}
\tau_{\ell-1}(\bs r),& \ell\le\ell',\\
\bs r,& \ell> \ell',
\end{cases}
\end{split}
\end{equation}
for target $\bs r$ with
\begin{equation}
 \tau_{\ell}(\bs r):=(x, \tau_{\ell}(y)).
\end{equation}
Then, the integral formulations for the reaction components $G_{\ell\ell'}^{\ast\star}$ can be written into the following uniform formulation
\begin{equation}
    G^{\ast\star}_{\ell\ell'}(\bs r, \bs r')=\frac{1}{2\pi }\int_{-\infty}^{\infty}\frac{\mathcal E_{\ell\ell'}(\lambda, \hat\tau_{\ell\ell'}^{\ast\star}(\bs r), \breve\tau_{\ell\ell'}^{\star}(\bs r'))}{2k_{\ell',y}}\sigma^{\ast\star}_{\ell \ell'}(\lambda) d\lambda,\quad \ast,\star=\uparrow,\downarrow,
\end{equation}
where
\begin{equation}\label{planewavekernel}
    \mathcal E_{\ell\ell'}(\lambda,\bs r, \bs r')= e^{\ri \lambda(x-x')+{\rm sgn}(y-y') (\ri k_{\ell,y}y-\ri k_{\ell',y}y')}.
\end{equation}

\begin{rem}
Note that the general reflection coefficients $\widetilde R_{\ell\ell'}(\lambda)$ in the top and bottom layers are defined as 
\begin{equation*}
	\widetilde R_{0,-1}(\lambda)=\widetilde R_{LL+1}(\lambda)\equiv 0.
\end{equation*}
Together with the formulas \eqref{disity1}-\eqref{disity2} imply that
\begin{equation*}
	\sigma_{0\ell^{\prime}}^{\downarrow\downarrow}(\lambda)=\sigma_{0\ell^{\prime}}^{\downarrow\uparrow}(\lambda)\equiv0,\quad 
	\sigma_{L\ell^{\prime}}^{\uparrow\uparrow}(\lambda)= \sigma_{L\ell^{\prime}}^{\uparrow\downarrow}(\lambda)\equiv0, \quad \ell'=0, 1, \cdots, L.
\end{equation*}
On the other hand, the definition of $\sigma_{\ell'\ell'}^{*\star}(\lambda)$ in \eqref{densityinsourcelayer} implies that
\begin{equation*}
	\sigma_{0 0}^{\uparrow\uparrow}(\lambda)=\sigma_{0 0}^{\downarrow\uparrow}(k_{\rho})=\sigma_{0 0}^{\downarrow\downarrow}(\lambda)\equiv 0,\quad \sigma_{L L}^{\uparrow\downarrow}(\lambda)=\sigma_{L L}^{\uparrow\uparrow}(k_{\rho})=\sigma_{\ell L}^{\downarrow\downarrow}(\lambda)\equiv 0,
\end{equation*}
which further leads to 
\begin{equation*}
	\sigma_{\ell 0}^{\uparrow\uparrow}(\lambda)=\sigma_{\ell 0}^{\downarrow\uparrow}(\lambda)=0,\quad
	\sigma_{\ell L}^{\uparrow\downarrow}(\lambda)=v_{\ell L}^{\downarrow\downarrow}(\lambda)=0,\quad \ell=0, 1, \cdots, L,
\end{equation*}
due to the relations in \eqref{disity1} and \eqref{disity1}. Thus, the following reaction components are zero
\begin{equation*}
	\begin{split}
		G_{0\ell^{\prime}}^{\downarrow\downarrow}(\bs r, \bs r')=G_{0\ell^{\prime}}^{\downarrow\uparrow}(\bs r, \bs r')\equiv 0,\quad 
		G_{L\ell^{\prime}}^{\uparrow\uparrow}(\bs r, \bs r')= G_{L\ell^{\prime}}^{\uparrow\downarrow}(\bs r, \bs r')\equiv 0,\\
		G_{\ell 0}^{\uparrow\uparrow}(\bs r, \bs r')=G_{\ell 0}^{\downarrow\uparrow}(\bs r, \bs r')\equiv 0,\quad
		G_{\ell L}^{\uparrow\downarrow}(\bs r, \bs r')=G_{\ell L}^{\downarrow\downarrow}(\bs r, \bs r')\equiv 0,
	\end{split}
\end{equation*}
for all $\ell, \ell'=0, 1, \cdots, L$. Therefore, the total number of non-zero reaction components is $4L^2$.
\end{rem}

\subsection{Effective transmission distance of reaction fields}\label{sec_transmission_distance}
\begin{figure}[t]  
\centering
\subfigure[transmission distance for $\widehat{\bs G}_{\ell'\ell'}^{\uparrow\downarrow}$ and $\widehat{\bs G}_{\ell'\ell'}^{\uparrow\uparrow}$]{\includegraphics[width=\linewidth]{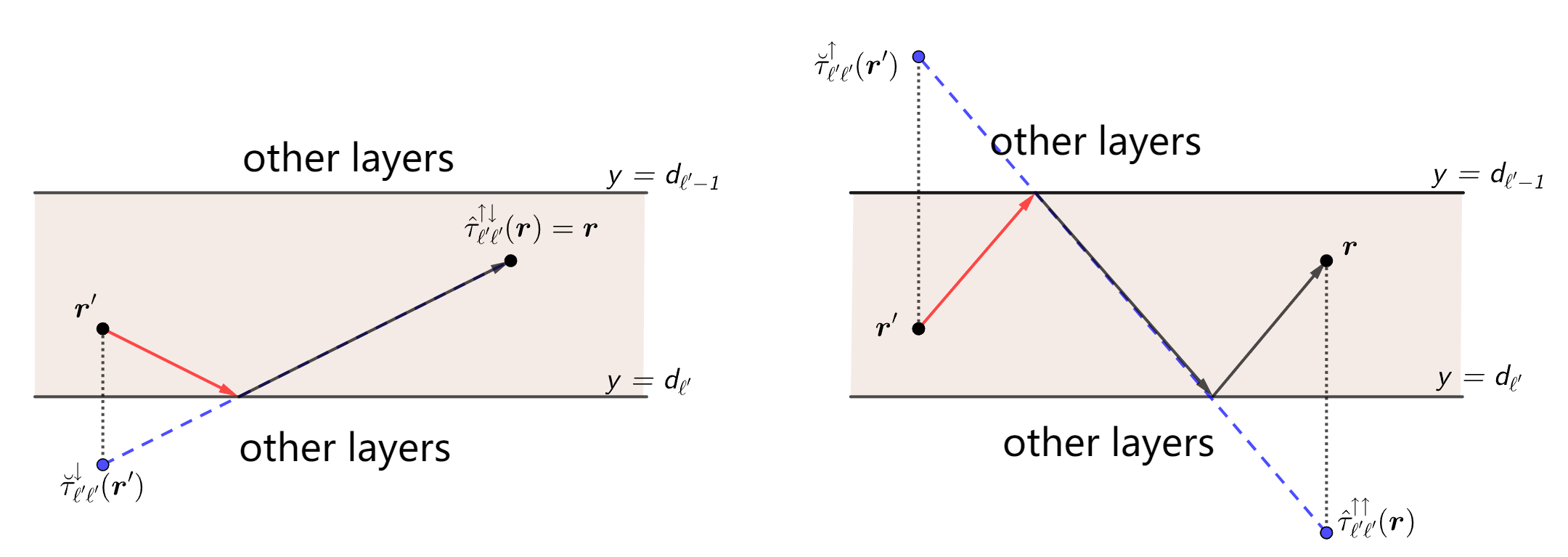}}
\subfigure[transmission distance for $\widehat{\bs G}_{\ell'\ell'}^{\downarrow\downarrow}$ and $\widehat{\bs G}_{\ell'\ell'}^{\downarrow\uparrow}$]{\includegraphics[width=\linewidth]{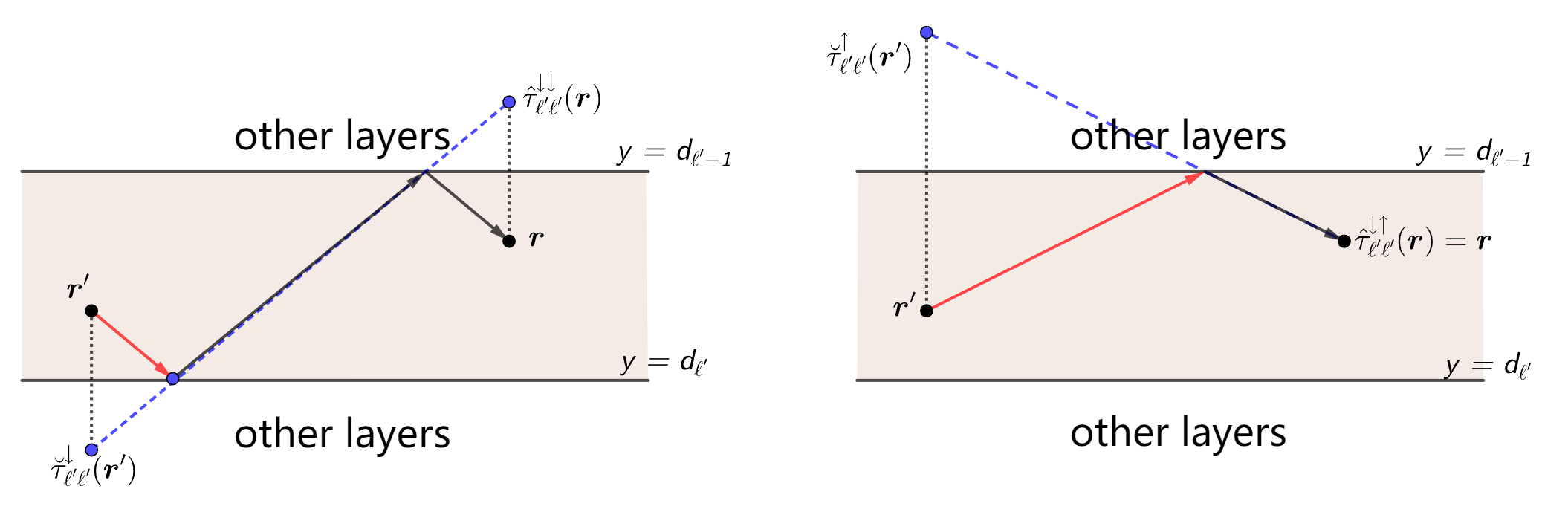}}
\caption{Equivalent polarization source coordinates and effective target locations in the case of $\ell=\ell'$.}
\label{targetoversource_eq}
\end{figure}

\begin{figure}[ht!]  
\centering
\subfigure[transmission distance for $\widehat{G}_{\ell\ell'}^{\uparrow\downarrow}$ and $\widehat{G}_{\ell\ell'}^{\uparrow\uparrow}$]{\includegraphics[width=\linewidth]{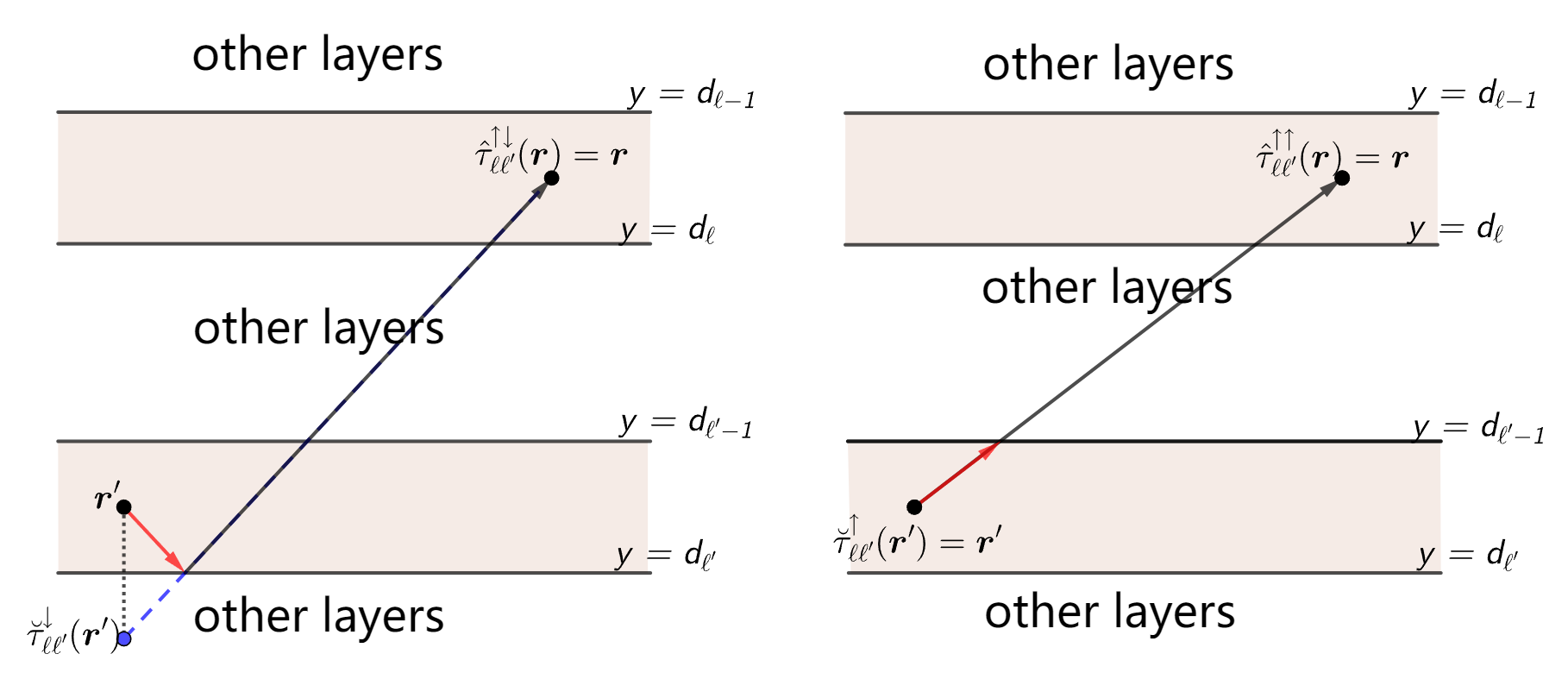}}
\subfigure[transmission distance for $\widehat{G}_{\ell\ell'}^{\downarrow\downarrow}$ and $\widehat{G}_{\ell\ell'}^{\downarrow\uparrow}$]{\includegraphics[width=\linewidth]{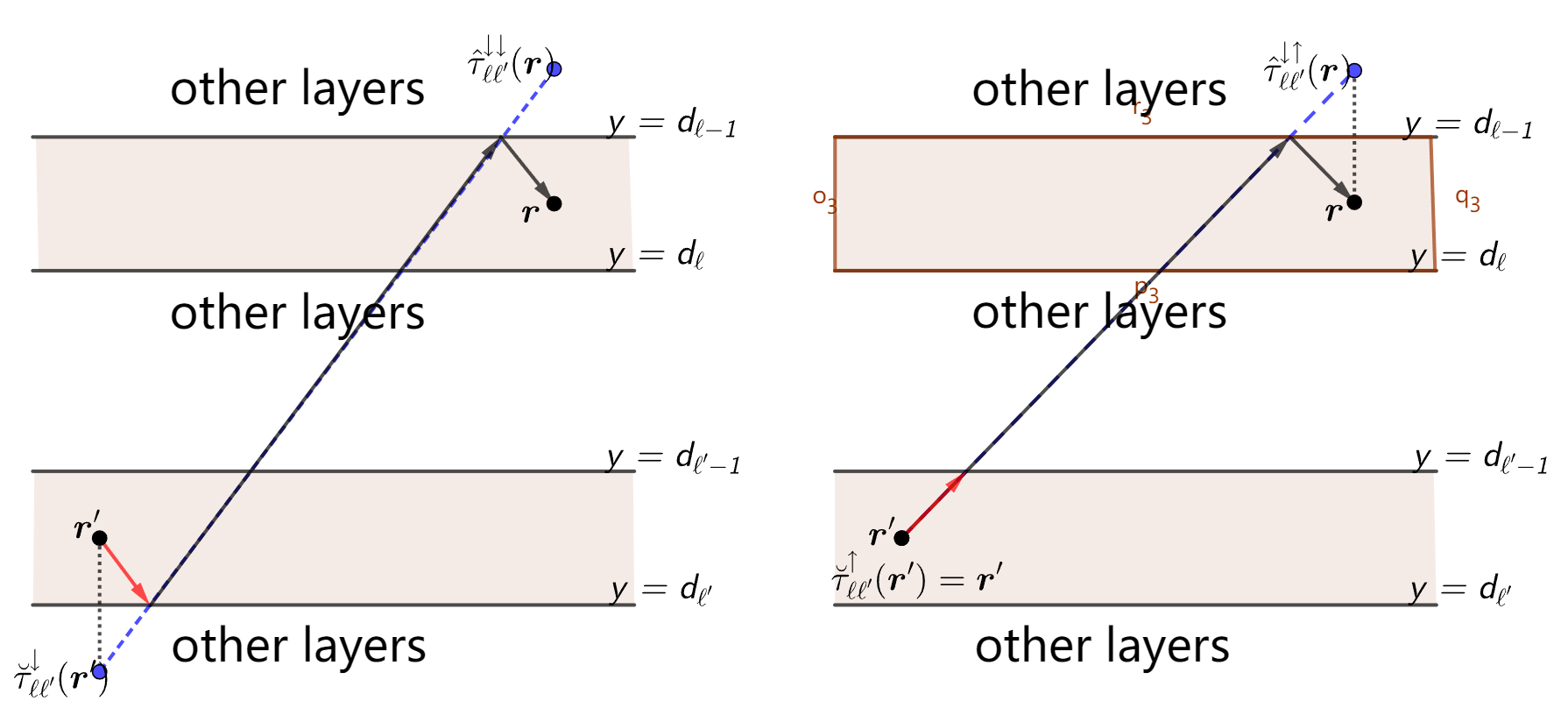}}
\caption{Equivalent polarization source coordinates and effective target locations in the case of $\ell<\ell'$.}
\label{fig_targetoversource}
\end{figure}

\begin{figure}[ht!]  
\centering
\subfigure[transmission distance for $\widehat{G}_{\ell\ell'}^{\uparrow\downarrow}$ and $\widehat{G}_{\ell\ell'}^{\uparrow\uparrow}$]{\includegraphics[width=\linewidth]{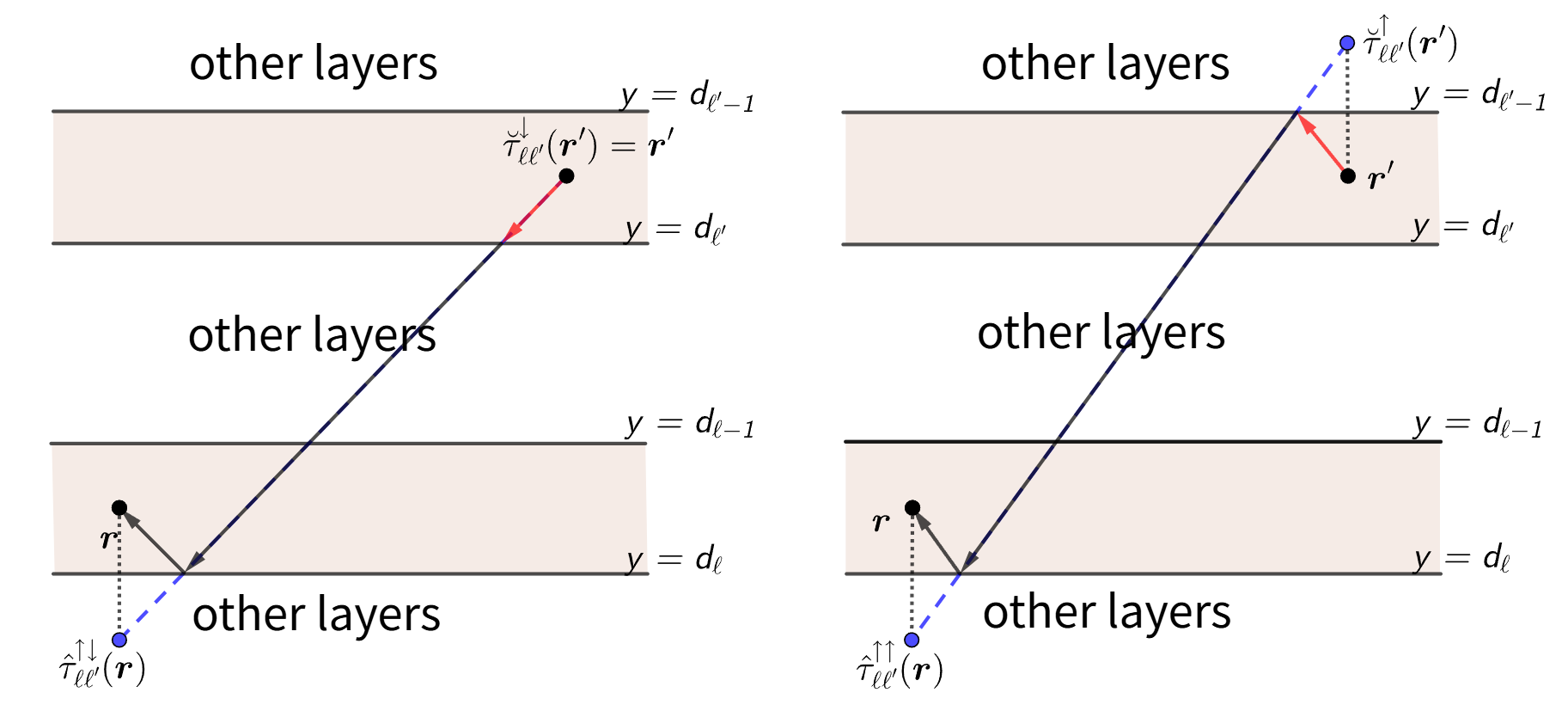}}
\subfigure[transmission distance for $\widehat{G}_{\ell\ell'}^{\downarrow\downarrow}$ and $\widehat{G}_{\ell\ell'}^{\downarrow\uparrow}$]{\includegraphics[width=\linewidth]{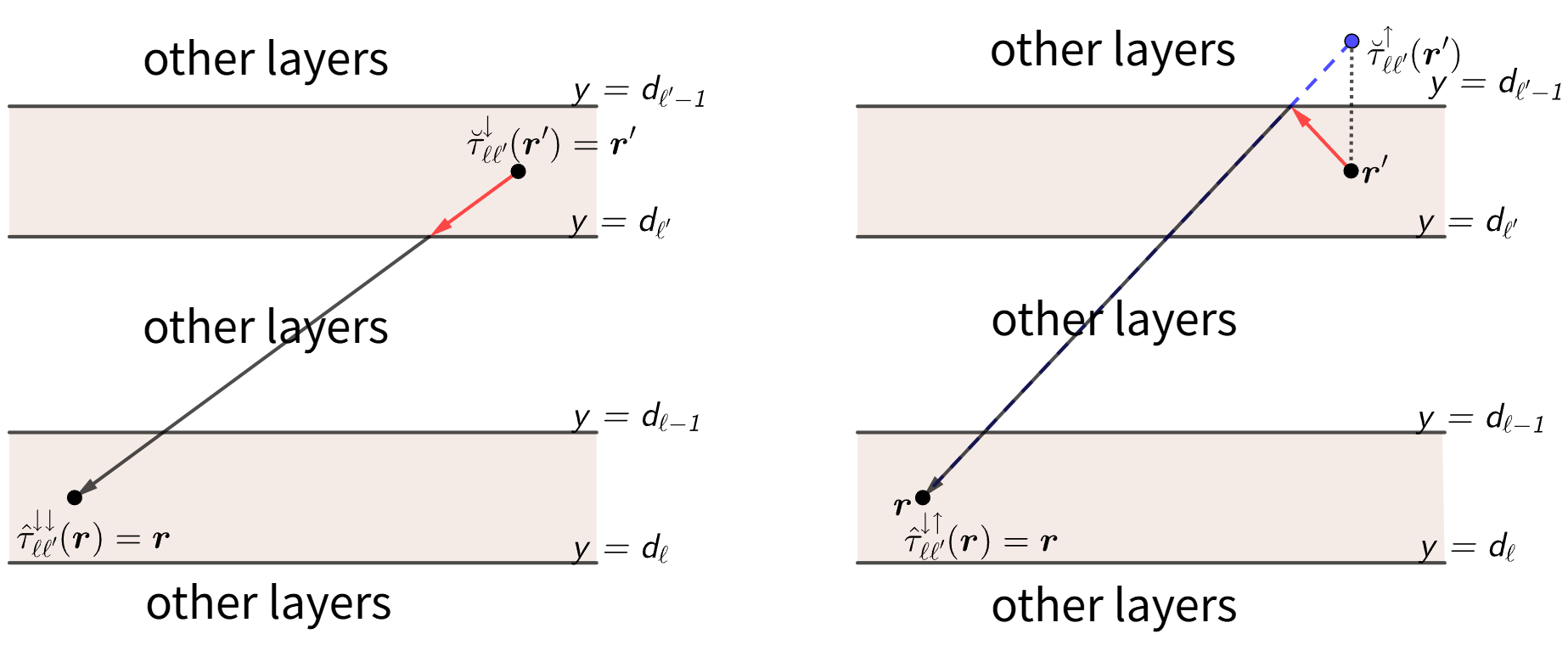}}
\caption{Equivalent polarization source coordinates and effective target locations in the case of $\ell>\ell'$.}
\label{fig_sourceovertarget}
\end{figure}

In addition to the equivalent polarization coordinates proposed in our previous works \cite{wang2019fast,zhang2020exponential,cai2025}, we introduce a new concept of effective location for the target particles, in order to account for the actual transmission distance of the reflected waves in layered media.
The upward and downward waves generated by the source at $\bs r'$ transmit to the target at $\bs r$ via different paths, see Figs. \ref{targetoversource_eq} to \ref{fig_sourceovertarget}, and generally induce upward and downward reaction fields. That is the physical background of the decomposition of four reaction field components.
The reaction field decomposition \eqref{layeredGreensfun} of Green's functions consists of the free space interaction provided $\ell = \ell'$, as well as the reaction field.
When categorized by the upward/downward field propagation directions, the reaction field is decomposed into (up to) four terms in \eqref{fourreactions}, each $\widehat{{G}}_{\ell\ell'}^{\ast\star}$ representing one type with upward ($\uparrow$), downward ($\downarrow$) or both directions, with the first symbol $\ast$ indicating the direction of wave arriving at the target, and the second symbol  $\star$  indicating the direction of the wave leaving the source.

When the target and the source come from the same layer, i.e. $\ell = \ell'$, waves of the reaction field must have at least one reflection on interfaces due to the subtraction of the free-space part, see Fig. \ref{targetoversource_eq} for an illustration. For instance, the reaction field component ${G}_{\ell'\ell'}^{\downarrow\downarrow}$ is interpreted as the superposition of waves that are downward at $\bs r'$ and downward at $\bs r$, including the wave marked by the solid line with two reflections in Fig. \ref{targetoversource_eq}(b), as well as any other contributions that may have experienced more reflections and transmissions on the interfaces. The minimal vertical transmission distance of these waves is given by the distance between the effective target location $\hat\tau_{\ell'\ell'}^{\downarrow\downarrow}(\bs r) = \tau_{\ell'-1}(\bs r)$ and the equivalent polarization source $\breve\tau_{\ell'\ell'}^{\downarrow}(\bs r') = \tau_{\ell'}(\bs r')$, as shown by the dashed line.
Indeed, in the exponent of \eqref{planewavekernel},
\begin{equation*}
\mathcal E_{\ell'\ell'}(\lambda, \hat\tau_{\ell'\ell'}^{\downarrow\downarrow}(\bs r), \breve\tau_{\ell'\ell'}^{\downarrow}(\bs r') ) \sim e^{-\lambda(\tau_{\ell'-1}(y) - \tau_{\ell'}(y'))}
\end{equation*}
vanishes exponentially as $\lambda\to+\infty$. 

The cases in which the target and the source are located in different layers are illustrated by Fig. \ref{fig_targetoversource} for $\ell < \ell'$, and Fig. \ref{fig_sourceovertarget} for $\ell > \ell'$, respectively.
The only difference with the previous case is that one of the reaction field components has minimal vertical transmission distance given by $|y - y'|$. This component must exist because it is \emph{not} equivalent to the free-space interaction due to the transmission through multiple layers.

In later discussion of the FMM implementation, the ``effective" field transmission distance
\begin{equation}\label{eq_transdict}
    d_{\ell\ell'}^{\ast\star}(\bs{r}, \bs{r}') = \left| \hat\tau_{\ell\ell'}^{\ast\star}(\bs r) - \breve\tau_{\ell\ell'}^{\star}(\bs r') \right|
\end{equation}
will be used as the criterion of far-field expansions.
Note that the field transmission distance is not shorter than that based on equivalent polarization source alone in our previous works \cite{wang2019fast,zhang2020exponential}, suggesting that wave sources in layered media are even more separated than previously thought.

\subsection{Indirect boundary integral equation and Nystr\"{o}m method}
By using the layered Green's function  $G_{\ell \ell'}(\bm r, \bm r')$ presented above, the solution of \eqref{extdirichleteq1} has boundary integral representation
\begin{equation}\label{integralrep}
	 u^s(\bs r')=\sum_{\ell=0}^{L}\frac{\eta_{\ell}}{\eta_{\ell'}}\int_{\Gamma_{\ell}}\phi(\bm {r})G_{\ell \ell'}(\bm r, \bm r')\,\mathrm{d}s,\quad \bs r'\in \Upsilon_{\ell'}\backslash{\Omega}_{\ell'}.
\end{equation}
It is worthy to point out that the integral is only on the boundaries of the scatterers, due to the usage of the layered Green's function. 
Applying boundary condition on $\{\Gamma_{\ell}\}_{\ell=0}^{L}$ gives integral equations
\begin{equation}\label{BoundaryIntergateEquation}
\sum_{\ell=0}^{L}\eta_{\ell}\int_{\Gamma_{\ell}}\phi(\bm {r})G_{\ell \ell'}(\bm r, \bm r')\,\mathrm{d}s=\eta_{\ell'}g(\bm {r}'),\quad \bs r'\in\Gamma_{\ell'},\; \ell'=0, 1, \cdots, L.
\end{equation}

For the discretization of the integral equations in \eqref{BoundaryIntergateEquation}, all boundaries $\{\Gamma_{\ell}\}_{\ell=0}^L$ are discretized into line segments as 
$$\Gamma_{\ell}\approx\bigcup_{e=1}^{N_{\ell}}\Gamma_{\ell e},\quad \ell=0, 1, \cdots, L,$$
where \(N_{\ell}\) denotes the number of elements in the mesh of \(\Gamma_{\ell}\). The total number of elements is denoted by $N = \sum_{\ell=0}^{L} N_{\ell}$. On each line element $\Gamma_{\ell e}$, we approximate the density function $\phi(\bs r)$ by a constant function $\phi_{\ell e}$. Substituting into \eqref{BoundaryIntergateEquation} and approximating the boundaries by $\{\Gamma_{\ell e}\}_{e=1}^{N_{\ell}}, \ell=0, 1, \cdots, L$, we obtain the linear system
\begin{equation}\label{linearSystem}    \displaystyle\sum_{\ell=0}^{L}\eta_{\ell}\sum_{e=1}^{N_{\ell}}\left(\int_{\Gamma_{\ell e}}G_{\ell \ell'}(\bm r, \bm c_{\ell' i})\,\mathrm{d}s\right)\phi_{\ell e} =\eta_{\ell'} g(\bm{c}_{\ell' i}),  \quad \bm{c}_{\ell' i}\in \Gamma_{\ell' i},      
\end{equation}
for $i=1,2,\cdots,N_{\ell'}$, $\ell' = 0,1,\cdots,L$, where the collocation point $\bm c_{\ell' i}$ is set to be the center of the element $\Gamma_{\ell' i}$.
It can be written as matrix formulation
\begin{equation}\label{EquationSet}
    \mathbb K \boldsymbol{\Phi}=\mathbf{b},
\end{equation}
where
\begin{equation}\label{matrixandvec}
    \mathbb{K}=\begin{bmatrix}
    \mathbb K_{00} & \mathbb K_{01} & \cdots & \mathbb K_{0L}\\
    \mathbb K_{10} & \mathbb K_{11} & \cdots & \mathbb K_{1L}\\
    \vdots & \vdots & \ddots & \vdots\\
    \mathbb K_{L0} & \mathbb K_{L1} & \cdots & \mathbb K_{LL}
\end{bmatrix},\quad \boldsymbol{\Phi}=\begin{bmatrix}
    \boldsymbol{\Phi}_{0} \\
    \boldsymbol{\Phi}_{1} \\
    \cdots \\
    \boldsymbol{\Phi}_{L}
\end{bmatrix},\quad
\mathbf{b}=\begin{bmatrix}
    \mathbf{b}_0\\
    \mathbf{b}_1\\
    \cdots \\
    \mathbf{b}_{L}
\end{bmatrix},
\end{equation}
with dense blocks 
\begin{equation}\label{matrixelmentsexpression}
    \mathbb{K}_{\ell\ell'}=(K_{ei}^{\ell\ell'})_{N_{\ell'}\times N_{\ell}},\quad K_{ei}^{\ell\ell'} = \int_{\Gamma_{\ell e}}G_{\ell \ell'}(\bm r, \bm c_{\ell' i})\,\mathrm{d}s,
\end{equation}
and
$$\boldsymbol{\Phi}^{\rm T}_{\ell}=\begin{bmatrix}
    \eta_{\ell}\phi_{\ell 1} & \eta_{\ell}\phi_{\ell 2} & \cdots &\eta_{\ell}\phi_{\ell,N_{\ell}}
\end{bmatrix},\quad \boldsymbol{b}^{\rm T}_{\ell}=\begin{bmatrix}
    \eta_{\ell}g(c_{\ell 1}) &  \eta_{\ell}g(c_{\ell 2}) & \cdots & \eta_{\ell}g(c_{\ell,N_{\ell}})
\end{bmatrix},$$
for $\ell,\ell'=0, 1, \cdots, L$.

Substituting the layered Green's function (\ref{layeredGreensfun}) into the expression(\ref{matrixelmentsexpression}) for \(K_{ei}^{\ell\ell'}\), gives
\begin{equation}\label{singleintegrate}
    K_{ei}^{\ell\ell'} = \delta_{\ell\ell'}\int_{\Gamma_{\ell e}}G^f_{\ell}(\bm r,\bm c_{\ell'i})\,\mathrm{d}s+\int_{\Gamma_{\ell e}}G^{\rm r}_{\ell\ell'}(\boldsymbol{r},\bm c_{\ell' i})\,\mathrm{d}s.
\end{equation}
Let us first consider the computation of the diagonal entries $K_{ii}^{\ell'\ell'}, \ell'=0, 1, \cdots, L; i=1, 2, \cdots, N_{\ell'}$. Note that $G^f_{\ell'}(\bm r,\bm c_{\ell' i})$ has weakly singularity at $\bs r=\bm c_{\ell' i}$ while $G^{\rm r}_{\ell'\ell'}(\boldsymbol{r},\bm c_{\ell' i})$ is smooth in the element $\Gamma_{\ell' i}$.
Therefore, the second term in the equation (\ref{singleintegrate}) can be directly approximated by 
\begin{equation}
    \int_{\Gamma_{\ell' i}}G^{\rm r}_{\ell'\ell'}(\boldsymbol{r},\bm c_{\ell' i})\,\mathrm{d}s\approx G^{\rm r}_{\ell'\ell'}(\bm c_{\ell'i},\bm c_{\ell' i})|\Gamma_{\ell' i}|.
\end{equation}
For weakly singular term, utilizing the asymptotic behavior
\[H_0^{(1)}(z) \sim 1 + \frac{2i}{\pi} \left( \ln\frac{z}{2} + \gamma\right), \quad z \to 0,\]
where \(\gamma\) denotes the Euler constant,
we decompose the integral into weakly singular and regular terms as follows
\begin{equation*}
\begin{aligned}
    \int_{\Gamma_{\ell' i}}G^f_{\ell'}(\bm r,\bm c_{\ell' i})\,\mathrm{d}s =&\frac{\ri}{4}\int_{\Gamma_{\ell' i}}\left[H_0^{(1)}(k_{\ell'}|\bs r-\bm c_{\ell' i}|) - 1 - \frac{2\ri}{\pi} \left( \ln\frac{k_{\ell'}|\bs r-\bm c_{\ell'i}|}{2} + \gamma\right)\right]\,\mathrm{d}s \\
    +&\frac{\ri}{4}\int_{\Gamma_{\ell'i}} \left[1 + \frac{2\ri}{\pi} \left( \ln\frac{k_{\ell'}|\bs r-\bm c_{\ell'i}|}{2} + \gamma\right)\right]\,\mathrm{d}s\\
    :=&I_{\ell'i}^r(\boldsymbol r)+I_{\ell'i}^s(\boldsymbol r).
\end{aligned}
\end{equation*}
For the regular term \(I_{\ell' i}^r(\bm r)\), the rectangle method yields approximation \(I_{\ell' i}^r(\bm c_{\ell e}) \approx 0\). For the weakly singular term \(I_{\ell'i}^s(\boldsymbol r)\), it is analytically tractable through integration by parts, i.e.,
\begin{equation}\label{ApproximateComputation}
    I_{\ell' i}^s(\boldsymbol c_{\ell' i}) = \left[\frac{\ri}{4}-\frac{1}{2\pi}\left(\gamma + \ln{\frac{k_{\ell'}|\Gamma_{\ell'i}|}{4}}-1\right)\right] |\Gamma_{\ell' i}|.
\end{equation}

For all off-diagonal entries, i.e., \((\ell, e)\neq (\ell', i)\), the integrand is regular in the integral domain and we are able to use the rectangle method for approximation. In summary, we shall use approximations
\begin{equation}\label{matrixelmemts}
 K_{ei}^{\ell\ell'} \approx \begin{cases}
 \displaystyle G_{\ell' \ell'}^{r}(\bm c_{\ell' i}, \bm c_{\ell' i})|\Gamma_{\ell' i}| +I_{\ell' i}^s(\boldsymbol c_{\ell' i}),\quad (\ell, e)= (\ell', i),\\
     \displaystyle G_{\ell \ell'}(\bm c_{\ell e}, \bm c_{\ell' i})|\Gamma_{\ell e}|, \quad (\ell, e)\neq (\ell', i).
 \end{cases}
 \end{equation}

 According to \eqref{matrixandvec}, \eqref{matrixelmentsexpression} and \eqref{matrixelmemts}, the coefficient matrix \(\mathbb{K}\) is a non-symmetric dense matrix. Consequently, the linear system (\ref{EquationSet}) can be solved by direct methods such as LU factorization or Gaussian elimination, though these require \(\mathcal{O}(N^3)\) operations. Alternatively, Krylov subspace iterative methods including GMRES, BiCG, could be employed. Therefore, the core computational task involves computing the matrix-vector product \(\boldsymbol{\Phi}^{(k+1)} = \mathbb{K}\boldsymbol{\Phi}^{(k)}\) at each iteration \(k\), where the iterate vector is defined as
 \[\boldsymbol{\Phi}^{(k)} = \begin{bmatrix}
     \boldsymbol{\Phi}_0^{(k)} & \boldsymbol{\Phi}_1^{(k)}& \cdots & \boldsymbol{\Phi}^{(k)}_L
 \end{bmatrix}^{\rm T}.\]
 Denoted the resultant vector by \(\boldsymbol{\Phi}^{(k+1)}=(\Phi_{\ell'i}^{(k+1)})_{N \times 1}\), each component \(\Phi_{\ell'i}^{(k+1)}\) is computed as follows
 \begin{equation}\label{Llayerscomponent}
    \begin{split}
        \Phi_{\ell'i}^{(k+1)}
        \approx \sum_{\ell=0}^{L}\sum\limits_{e=1}^{N_{\ell}}G^r_{\ell \ell'}(\bm c_{\ell e},\bm c_{\ell' i})Q_{\ell e}+\sum\limits_{\substack{e=1\\e\neq i}}^{N_{\ell'}}G^f_{\ell'}(\bm c_{\ell e},\bm c_{\ell' i})Q_{\ell e} 
        +I_{\ell' i}^s(\boldsymbol c_{\ell' i})\eta_{\ell'}\phi_{\ell' i}^{(k)},
    \end{split}
\end{equation}
where $Q_{\ell e}=\eta_{\ell}\phi_{\ell e}^{(k)}|\Gamma_{\ell e}|$. Here, the approximations \(\eqref{matrixelmemts}\) and (\eqref{ApproximateComputation}) for \(K_{ie}^{\ell\ell'}\) have been used. Define
\begin{equation}\label{componentwise}
\Phi_{\ell'i}^{f}=\sum\limits_{\substack{e=1\\e\neq i}}^{N_{\ell'}}G^f_{\ell'}(\bm c_{\ell e},\bm c_{\ell' i})Q_{\ell e},\quad \Phi_{\ell'i}^{r}=\sum_{\ell=0}^{L}\sum\limits_{e=1}^{N_{\ell}}Q_{\ell e}G^r_{\ell \ell'}(\bm c_{\ell e},\bm c_{\ell' i}).
\end{equation}
Then, the matrix product vector \eqref{Llayerscomponent} can be re-expressed as
\begin{equation}\label{decomposedpotential}
\Phi_{\ell'i}^{(k+1)}=\Phi_{\ell'i}^{f}+\Phi_{\ell'i}^{r}+I_{\ell'i}^s(\bm{c}_{\ell'i})\eta_{\ell'}\phi_{\ell' i}^{(k+1)}, \quad i= 1, 2,\cdots, N_{\ell'};\; \ell'=0, 1, \cdots, L.
\end{equation}
Direct computation incurs \(\mathcal{O}(N^2)\) cost for each iteration. Therefore, both direct and iterative methods become computationally infeasible for large-scale problems.

\section{Fast multipole acceleration and preconditioning}\label{sect4}
In this section, we have improved a fast multipole method (FMM) from our previous work (cf.\cite{zhang2020exponential,wang2019fast}), which accelerates matrix-vector products $\boldsymbol{\Phi}^{(k+1)}=\mathbb{K}\boldsymbol{\Phi}^{(k)}$ within the iterative solver. This technique reduces the computational cost per iteration to $\mathcal{O}(N\log N)$.

The formulation \(\eqref{decomposedpotential}\) shows that \(\{\Phi_{\ell'i}^{f}\}_{i = 1}^{N_{\ell'}}\) and \(\{\Phi_{\ell'i}^{r}\}_{i = 1}^{N_{\ell'}}\) can be computed separately. For free-space parts \(\{\Phi_{\ell'i}^{f}\}_{i = 1}^{N_{\ell'}}\), the classical FMM  can reduce the cost to \(\mathcal{O}(N_{\ell'}\log N_{\ell'})\) per layer. The main challenge therefore lies in accelerating the computation of the reaction field components \(\{\Phi_{\ell'i}^{r}\}_{i = 1}^{N_{\ell'}}\). 

From the reaction kernel \(G_{\ell\ell'}^{r}(\mathbf{r}, \mathbf{r}')\) in \eqref{reactioncomponent}-\eqref{fourreactions}, \(\Phi_{\ell'i}^{r}\) admits the decomposition 
\begin{equation}\label{reactionfieldcomponents}
\Phi_{\ell'i}^{r}=\sum\limits_{\ell=0}^{L}\big[\Phi_{\ell\ell'}^{\uparrow\downarrow}(\bs r_{\ell' i})+\Phi_{\ell\ell'}^{\uparrow\uparrow}(\bs r_{\ell' i})+\Phi_{\ell\ell'}^{\downarrow\uparrow}(\bs r_{\ell' i})+\Phi_{\ell\ell'}^{\downarrow\downarrow}(\bs r_{\ell' i})\big],
\end{equation}
where
\begin{equation}\label{unfiormintegral}
	\Phi_{\ell\ell'}^{*\star}(\bs r_{\ell' i})=\sum\limits_{e=1}^{N_{\ell}}Q_{\ell e}G^{*\star}_{\ell\ell'}(\bs r_{\ell e}, \bs r_{\ell' i})=\sum\limits_{e=1}^{N_{\ell}}\frac{Q_{\ell e}}{2\pi }\int_{-\infty}^{\infty}\frac{\mathcal E_{\ell\ell'}(\lambda, \hat\tau_{\ell\ell'}^{\ast\star}(\bs r_{\ell e}), \breve\tau_{\ell\ell'}^{\star}(\bs r'_{\ell'i}))}{2k_{\ell',y}}\sigma^{\ast\star}_{\ell \ell'}(\lambda) d\lambda,
\end{equation}
for $*,\star=\uparrow,\downarrow$. This is a uniform integral representation for any reaction components indicated by quadruple $(\ell,\ell',\ast,\star)$. Hence, we generally need to 
compute the following summations 
\begin{equation}\label{summations}
\begin{split}
	\Phi^+(\bs r_{i}, \sigma)&=\sum\limits_{j=1}^{M_s}Q_j'\mathcal I^{k'k}_{00}(x'_j-x_i, y'_j, y_i, \sigma),\quad i= 1, \cdots, M_t,\\
	\Phi^-(\bs r_{i}, \sigma)&=\sum\limits_{j=1}^{M_s}Q_j''\mathcal I^{kk''}_{00}(x_j''-x_i, y_i, y''_j, \sigma),\quad i= 1, \cdots, M_t,
\end{split}
\end{equation}
where the target points $\{\bs r_i=(x_i, y_i)\}_{i=1}^{M_t}$ and the source points $\{\bs r_j'=(x_j', y_j')\}_{j=1}^{M_s}$ or $\{\bs r_j''=(x_j'', y_j'')\}_{j=1}^{M_s}$ should be obtained by the reflections $\hat\tau_{\ell\ell'}^{\ast\star}$, $\breve{\tau}_{\ell\ell'}^{\star}$ of the original target and source points $\bs r_{\ell e}, \bs r_{\ell' i}$ accordingly, $k, k'$ and $\sigma$ is the triple $(k_{\ell}, k_{\ell'}, \sigma_{\ell\ell'}^{\ast\star})$ determined according to quadruple $(\ell,\ell',\ast,\star)$ and $\mathcal I^{kk'}_{00}(x, y, y', \sigma)$ is defined in \eqref{generalintegral}. Due to the reflections, we always have the source and target points are separate by one of the interface in the computation of $\Phi^{\pm}(\bs r_{i}, \sigma)$. Moreover, the reflection always ensure that $y_j'>y_i$ and $y_i>y_j''$.

\subsection{Fast multipole method for reaction components}
The basic theory for the FMM of two-dimensional Helmholtz equation in layered media is reported in \cite{zhang2020exponential}, but without using the concept of effective locations for the target points. In this subsection, we apply this theory to the potentials \eqref{summations} which are defined using both effective locations for targets and equivalent polarization coordinates for sources and provide details for the implementation of the FMM.

The key ingredients in the derivation of the expansion theory are the following expansions obtained by the formulation of the generating function of the cylindrical Bessel functions of the first kind. We refer to \cite{zhang2020exponential} for detailed derivation.
\begin{proposition}\label{prop:Funk-Hecke-limit}
	Given four points $\bs r=(x, y), \bs r'=(x', y')$, $\bs r_c=(x_c, y_c)$ and $\bs r_c'=(x'_c, y_c')$ in $\mathbb R^2$ such that $|\bs r-\bs r_c|<|\bs r'-\bs r_c|$ and $|\bs r-\bs r_c'|>|\bs r'-\bs r_c'|$.  Suppose $y>y'$, $y>y_c'$ and $y_c>y'$. Denoted by
\begin{equation}
    \mathcal E^{kk'}(x, y, y')=e^{\ri \lambda x+\ri \sqrt{k^2-\lambda^2}y-\ri \sqrt{k'^2-\lambda^2}y'},
\end{equation}
then there holds the following expansions
\begin{equation}\label{reactfieldMELE}
    \begin{aligned}
        &\displaystyle\mathcal E^{kk'}(x-x', y, y') \displaystyle=\sum\limits_{n=-\infty}^{\infty}J_{n}(k\rho_{c})e^{\ri n\theta_{c}}\omega(\lambda,k)^n\mathcal E^{kk'}(x_c-x', y_c, y'),\\
        &\displaystyle\mathcal E^{kk'}(x-x', y, y')
        \displaystyle=\sum\limits_{n=-\infty}^{\infty}(-1)^nJ_{n}(k'\rho_{c}')e^{\ri n\theta_{c}'}\omega(\lambda,k')^{n}\mathcal E^{kk'}(x-x_c', y, y_c'),
    \end{aligned}
\end{equation}
where 
\begin{equation}
w(\lambda, k)=\frac{\sqrt{k^{2}-\lambda^{2}}+\ri\lambda}{k},
\end{equation}
$(\rho_{c}', \theta_{c}')$ and $(\rho_{c}, \theta_{c})$ are the polar coordinates of $\bs r'-\bs r'_c$ and $\bs r-\bs r_c$, respectively.
\end{proposition}

Define
\begin{equation}\label{generalintegral}
	\begin{split}
		\mathcal I^{kk'}_{nm}(x, y, y', \sigma):=\int_{-\infty}^{\infty}\mathcal E^{kk'}(x, y, y')\omega(\lambda,k)^{n}\omega(\lambda,k')^{m}\sigma(\lambda) d\lambda\quad n,m\in\mathbb{Z}.
	\end{split}
\end{equation}
Then using the expansions \eqref{reactfieldMELE} inside the integral and then exchange the order of the infinite summation and integral, we obtain
\begin{equation}\label{reactfieldLE1}
	\mathcal I^{kk'}_{0m}(x-x', y, y', \sigma)=\sum\limits_{n=-\infty}^{\infty}J_{n}(k\rho_{c})e^{\ri n\theta_{c}}\mathcal I_{nm}^{kk'}(x_c-x', y_c, y',\sigma),\quad m=0, 1, \cdots,
\end{equation}
and
\begin{equation}\label{reactfieldME1}
	\mathcal I^{kk'}_{m0}(x-x', y, y', \sigma)=\sum\limits_{n=-\infty}^{\infty}(-1)^nJ_{n}(k'\rho'_c)e^{\ri n\theta'_c}\mathcal I_{mn}^{kk'}(x-x_c', y, y_c',\sigma),\quad m=0, 1, \cdots,
\end{equation} 
respectively. The convergence of the above expansions has been proved and verified in \cite{zhang2020exponential}.

Given source center $(x'_c, y'_c)$, $(x''_c, y''_c)$ close to the source points, the expansion \eqref{reactfieldLE1} and \eqref{reactfieldME1} gives
\begin{equation}\label{mesingle}
\begin{split}
\mathcal I^{k'k}_{00}(x'_j-x_i, y'_j, y_i, \sigma)=&\sum\limits_{n=-\infty}^{\infty}J_{n}(k'\rho_{j}')e^{\ri n\theta_{j}'}\mathcal I_{n0}^{k'k}(x_c'-x_i, y_c', y_i,\sigma),\\
\mathcal I^{kk''}_{00}(x_j''-x_i, y_i, y''_j, \sigma)
=&\sum\limits_{n=-\infty}^{\infty}(-1)^nJ_{n}(k''\rho_{j}'')e^{\ri n(\pi-\theta_{j}'')}\mathcal I_{0n}^{kk''}(x_c''-x_i,y_i, y_c'', \sigma),\\
=&\sum\limits_{n=-\infty}^{\infty}J_{n}(k''\rho_{j}'')e^{-\ri n\theta_{j}''}\mathcal I_{0n}^{kk''}(x_c''-x_i,y_i, y_c'', \sigma),
\end{split}
\end{equation}
where $(\rho_{j}', \theta_{j}')$ and $(\rho_{j}'', \theta_{j}'')$ are the polar coordinates of $(x_j'-x_c', y_j'-y_c')$ and $(x_j''-x_c'', y_j''-y_c'')$, respectively. For the derivation of the second expansion, we actually use \eqref{reactfieldME1} with $(x, y)=(-x_i, y_i)$, $(x', y')=(-x_j'', y_j'')$ and $(x_c', y_c')=(-x_c'', y_c'')$.

Therefore, we obtain multipole expansions (ME):
\begin{equation}\label{MEformula}
\begin{split}
\Phi^+(\bs r_{i}, \sigma)=&\sum\limits_{|n|=0}^{\infty}\alpha_{n}^{+}\mathcal I_{n0}^{k'k}(x_c'-x_i, y_c', y_i,\sigma),\\
\Phi^-(\bs r_{i}, \sigma)=&\sum\limits_{|n|=0}^{\infty}\alpha_{n}^{-}\mathcal I_{0,-n}^{kk''}(x_c''-x_i,y_i, y_c'', \sigma),
\end{split}
\end{equation}
where
\begin{equation}\label{MEcoefficients}
 \alpha_{n}^{+}=\sum\limits_{j=1}^{M_s}Q'_jJ_{n}(k'\rho_{j}')e^{\ri n\theta_{j}'},\quad
 \alpha_{n}^{-}=\sum\limits_{j=1}^{M_s}(-1)^nQ''_jJ_{n}(k''\rho_{j}'')e^{\ri n\theta_{j}''}.
\end{equation}

Similarly, given the target center $(x_c, y_c)$ close to the target points $\{(x_i, y_i)\}$, the  expansion \eqref{reactfieldLE1} and \eqref{reactfieldME1} gives the local expansions (LE):
\begin{equation}\label{LEformula}
\Phi^+(\bs r_{i}, \sigma)=\sum\limits_{|n|=0}^{\infty}\beta_{n}^+J_{n}(k\rho_{i})e^{\ri n\theta_{i}},\quad
\Phi^-(\bs r_{i},\sigma)=\sum\limits_{|n|=0}^{\infty}\beta_{n}^-J_{n}(k\rho_{i})e^{\ri n\theta_{i}},
\end{equation}
where $(\rho_{i},\theta_{i})$ is the polar coordinates of $(x_i-x_c, y_i-y_c)$, and
\begin{equation}
\beta_{n}^+=\sum\limits_{j=1}^{M_s}Q'_j(-1)^n\mathcal I_{0n}^{k'k}(x_j'-x_c, y_j', y_c,\sigma),\quad
\beta_{n}^-=\sum\limits_{j=1}^{M_s}Q''_j\mathcal I_{-n0}^{kk''}(x_j''-x_c,y_c, y_j'', \sigma).
\end{equation}

Let $(\tilde x_c', \tilde y_c')$ be another source center close to $(x_c', y_c')$ and $(\tilde x_c'', \tilde y_c'')$ be another source center close to $(x_c'', y_c'')$. Then, the multipole expansions with respect to new source centers $(\tilde x_c', \tilde y_c')$ and $(\tilde x_c'', \tilde y_c'')$ are given by
\begin{equation}
\begin{split}
\Phi^+(\bs r_{i}, \sigma)=&\sum\limits_{|n|=0}^{\infty}\tilde \alpha_{n}^+\mathcal I_{n0}^{kk'}(\tilde x_c'-x_i, \tilde y_c', y_i,\sigma),\\
\Phi^-(\bs r_{i}, \sigma)=&\sum\limits_{|n|=0}^{\infty}\tilde\alpha_{n}^- \mathcal I_{0,{-n}}^{kk''}(\tilde x_c''-x_i,y_i, \tilde y_c'', \sigma),
\end{split}
\end{equation}
where
\begin{equation}
\begin{split}
 \tilde\alpha_{n}^+=\sum\limits_{j=1}^{M_s}Q'_jJ_{n}(k'\tilde\rho_{j}')e^{\ri n\tilde\theta_{j}'},\quad 
 \tilde\alpha_{n}^-=\sum\limits_{j=1}^{M_s}(-1)^nQ''_jJ_{n}(k''\tilde\rho_{j}'')e^{\ri n\tilde\theta_{j}''},
\end{split}
\end{equation}
$(\tilde\rho'_{j},\tilde\theta'_{j})$ and $(\tilde\rho''_{j},\tilde\theta''_{j})$ are the polar coordinates of $(x'_j-\tilde{x}'_{c}, y'_j-\tilde{y}'_{c})$ and $(x''_j-\tilde{x}''_{c}, y''_j-\tilde{y}''_{c})$.
By Graf's addition theorem (cf. \cite{martin2006multiple}), we have M2M shifting operators:
\begin{equation}\label{metome}	\tilde\alpha_{n}^+=\sum_{|m|=0}^{\infty}J_{n-m}(k\rho'_{c})e^{\ri(n-m)\theta'_{c}}{\alpha}_{m}^+,\quad
\tilde\alpha_{n}^-=\sum_{|m|=0}^{\infty}J_{m-n}(k''\rho''_{c})e^{\ri(n-m)\theta''_{c}}{\alpha}_{m}^-,
\end{equation}
where $(\rho'_{c},\theta'_{c})$ and $(\rho''_{c},\theta''_{c})$ are the polar coordinates of $(x'_c-\tilde{x}'_c, y_c-\tilde{y}'_c)$ and  $({x}''_c-\tilde{x}''_c, {y}''_c-\tilde{y}''_c)$. 

Let $(\tilde x_c, \tilde y_c)$ be another target center close to $(x_c, y_c)$. By the local expansions in \eqref{LEformula} and Graf's addition theorem (cf. \cite{martin2006multiple}), we have
\begin{equation}
	\Phi^+(\bs r_{i}, \sigma)=\sum\limits_{|n|=0}^{\infty}\tilde{\beta}_{n}^+J_{n}(k\tilde{\rho}_{i})e^{\ri n\tilde\theta_i},\quad \Phi^-(\bs r_{i}, \sigma)=\sum\limits_{|n|=0}^{\infty}\tilde\beta_{n}^-J_{n}(k\tilde\rho_{i})e^{\ri n\tilde\theta_i},
\end{equation}
where 
\begin{equation}\label{letole1}
	\tilde{\beta}_{n}^+=\sum_{|m|=0}^{\infty}J_{m}(k\tilde\rho_{c})e^{\ri m\tilde\theta_{c}}\beta_{n+m}^+,\quad \tilde\beta_n^-=\sum_{|m|=0}^{\infty}(-1)^mJ_{m}(k\tilde\rho_{c})e^{-\ri m\tilde\theta_{c}}\beta_{m-n}^-,
\end{equation}
are the L2L shifting operators and $(\tilde\rho_{i},\tilde\theta_{i})$ and $(\tilde\rho_{c},\tilde\theta_{c})$ are the polar coordinates of $(x_i-\tilde{x}_c, y_i-\tilde{y}_c)$ and $(x_c-\tilde{x}_c, y_c-\tilde{y}_c)$.

Next, we consider the multipole-to-local translation. Applying expansion formulations \eqref{reactfieldLE1}-\eqref{reactfieldME1} to the multipole expansions in \eqref{mesingle}, we obtain
\begin{equation*}
	\begin{split}
		\psi^+(\bs r_{i}, \sigma)=&\sum\limits_{|m|=0}^{\infty}\sum\limits_{|n|=0}^{\infty}\alpha_{n}^+(-1)^mJ_m(k\rho_i)e^{\ri m\theta_i}\mathcal I_{nm}^{k'k}( x'_c-x_c,y'_c,y_c,\sigma),\\
		\psi^-(\bs r_{i}, \sigma)=&\sum\limits_{|m|=0}^{\infty}\sum\limits_{|n|=0}^{\infty}\alpha_{n}^-(-1)^mJ_m(k\rho_i)e^{-\ri m\theta_i}\mathcal I_{mn}^{kk''}( x''_c-x_c,y_c,y''_c,\sigma).
	\end{split}
\end{equation*}
Matching with the local expansions in \eqref{LEformula} gives M2L translations
\begin{equation}\label{metoleimage}
	\beta_m^+=\sum\limits_{|n|=0}^{\infty}(-1)^m\mathcal I_{nm}^{k'k}( x'_c-x_c,y'_c,y_c,\sigma)\alpha_{n}^+,\quad
	\beta_m^-=\sum\limits_{|n|=0}^{\infty}\mathcal I_{-mn}^{kk''}( x''_c-x_c,y_c,y''_c,\sigma)\alpha_{n}^-.
\end{equation}

Using the truncated expansions, shifting and translation operators in the framework of the classic FMM, we implement a FMM for fast calculation of $\{\Phi^{\pm}(\bs r_i,\sigma)\}_{i=1}^{M_t}$ at any desired accuracy. The pseudo-code of the algorithm and the overall FMM for the computation of the interactions defined in \eqref{unfiormintegral} is similar to the algorithms we have reported in our previous work \cite{wang2019fast}. Our implementation actually used the adaptive tree structure by changing the FMM framework to its adaptive version (cf. \cite{Ying2004}).

\subsection{An overlapping domain decomposition preconditioner}
Obviously, equation (\ref{BoundaryIntergateEquation}) is a first-kind Fredholm integral equation. By the compactness of the integral operator and the spectral theory of compact operators, it has an infinite condition number. Although the discretized linear system (\ref{linearSystem}) is an approximation of the integral equation (\ref{BoundaryIntergateEquation}) and thus exhibits a finite condition number, its value grows unbounded as the degree of freedom increases. Such an ill-conditioning severely degrades the convergence rate of iterative solvers. To enable efficient large-scale simulations, we will use a preconditioning technique for capacitance extraction of conductors (cf. \cite{Nabors1991Fast}) to develop a tailored preconditioning framework for layered media problems. Specifically, high performance can be achieved by applying overlapping domain decomposition preconditioning within each layer.

Recall the approximations in \eqref{matrixelmemts}, each diagonal block of the matrix $\mathbb K$ has the form
\begin{equation*}
	\mathbb{K}_{\ell\ell}= \scalebox{0.75}{
	$	\begin{bmatrix}
			|\Gamma_{\ell 1}|G_{\ell \ell}^r(\bs c_{\ell1}, \bs c_{\ell1})+I_{\ell1}^s(\bm{c}_{\ell1}) & |\Gamma_{\ell2}|G_{\ell \ell}(\bs c_{\ell2}, \bs c_{\ell1}) & \cdots &|\Gamma_{\ell N_{\ell}}|G_{\ell \ell}(\bs c_{\ell N_{\ell}}, \bs c_{\ell1})\\
			|\Gamma_{\ell1}|G_{\ell \ell}(\bs c_{\ell1}, \bs c_{\ell2}) & |\Gamma_{\ell2}|G_{\ell \ell}^r(\bs c_{\ell2}, \bs c_{\ell2})+I_{\ell2}^s(\bm{c}_{\ell2}) & \cdots &|\Gamma_{\ell N_{\ell}}|G_{\ell \ell}(\bs c_{\ell N_{\ell}}, \bs c_{\ell2})\\
			\vdots & \vdots & \ddots & \vdots\\
			|\Gamma_{\ell1}|G_{\ell \ell}(\bs c_{\ell1}, \bs c_{\ell N_{\ell}})& |\Gamma_{\ell2}|G_{\ell \ell}(\bs c_{\ell2}, \bs c_{\ell N_{\ell}}) & \cdots &|\Gamma_{\ell N_{\ell}}|G_{\ell\ell}^{r}(\bs c_{\ell N_{\ell}}, \bs c_{\ell N_{\ell}})+I_{\ell N_{\ell}}^s(\bm{c}_{\ell N_{\ell}})
		\end{bmatrix}$},
\end{equation*}
for all $\ell=0, 1, \cdots, L$.  We ignore the reaction filed component \(G_{\ell \ell'}^r\) in the diagonal blocks $\mathbb K_{\ell\ell}$ and define
$$
\widetilde{\mathbb K}_{\ell\ell}=\begin{bmatrix}
    I_{\ell 1}^s(\bm{c}_{\ell 1}) & |\Gamma_{\ell2}|G^f_{\ell}(\bs c_{\ell2}, \bs c_{\ell1}) & \cdots &|\Gamma_{\ell N_{\ell}}|G^f_{\ell}(\bs c_{\ell N_{\ell}}, \bs c_{\ell1})\\
    |\Gamma_{\ell1}|G^f_{\ell}(\bs c_{\ell1}, \bs c_{\ell2})& I_{\ell 2}^s(\bm{c}_{\ell 2})   & \cdots &|\Gamma_{\ell N_{\ell}}|G^f_{\ell}(\bs c_{\ell N_{\ell}}, \bs c_{\ell2})\\
    \vdots & \vdots & \ddots & \vdots\\
    |\Gamma_{\ell1}|G^f_{\ell}(\bs c_{\ell1}, \bs c_{\ell N_{\ell}})& |\Gamma_{\ell2}|G^f_{\ell}(\bs c_{\ell2}, \bs c_{\ell N_{\ell}}) & \cdots &I_{\ell N_{\ell}}^s(\bm{c}_{\ell N_{\ell}})
\end{bmatrix},
$$
then, we obtain an approximation for diagonal part of $\mathbb K$, i.e.,
$$
\widetilde{\mathbb K}=\begin{bmatrix}
    \tilde{\mathbb K}_{00} & 0 & \cdots & 0\\
    0 & \tilde{\mathbb K}_{11} & \cdots & 0\\
    \vdots & \vdots & \ddots & \vdots\\
    0 & 0 & \cdots   & \tilde{\mathbb K}_{LL}
\end{bmatrix}.
$$
Although $\widetilde{\mathbb{K}}$ is block-diagonal, each block corresponds to an approximation of an free space scattering problem within a single layer, whose solution remains computationally expensive. To address this, we leverage the tree structure generated by the FMM  to construct an approximation of $\widetilde{\mathbb{K}}^{-1}$ as a preconditioner for solving the linear system \eqref{EquationSet}.

The spatial decay property of the free-space Green's function \eqref{freegreensfun} implies that $|G^f_{\ell}(\mathbf{r},\mathbf{r}')| \to 0$ as $|\mathbf{r} - \mathbf{r}'| \to \infty$. Consequently, matrix entries $|\Gamma_{\ell j}|G^f_{\ell}(\mathbf{c}_{\ell j}, \mathbf{c}_{\ell k})$ in $\widetilde{\mathbb{K}}_{\ell \ell}$ exhibit similar decay. This indicates that for any unknown $\phi_{\ell e}$, only strongly related to its nearby unknowns. Thus, in the construction of a preconditioner, influences from distant DOFs can be neglected, a proximity criterion naturally provided by the FMM tree structure.
\begin{figure}[ht!]  
\centering
\includegraphics[width=0.5\linewidth]{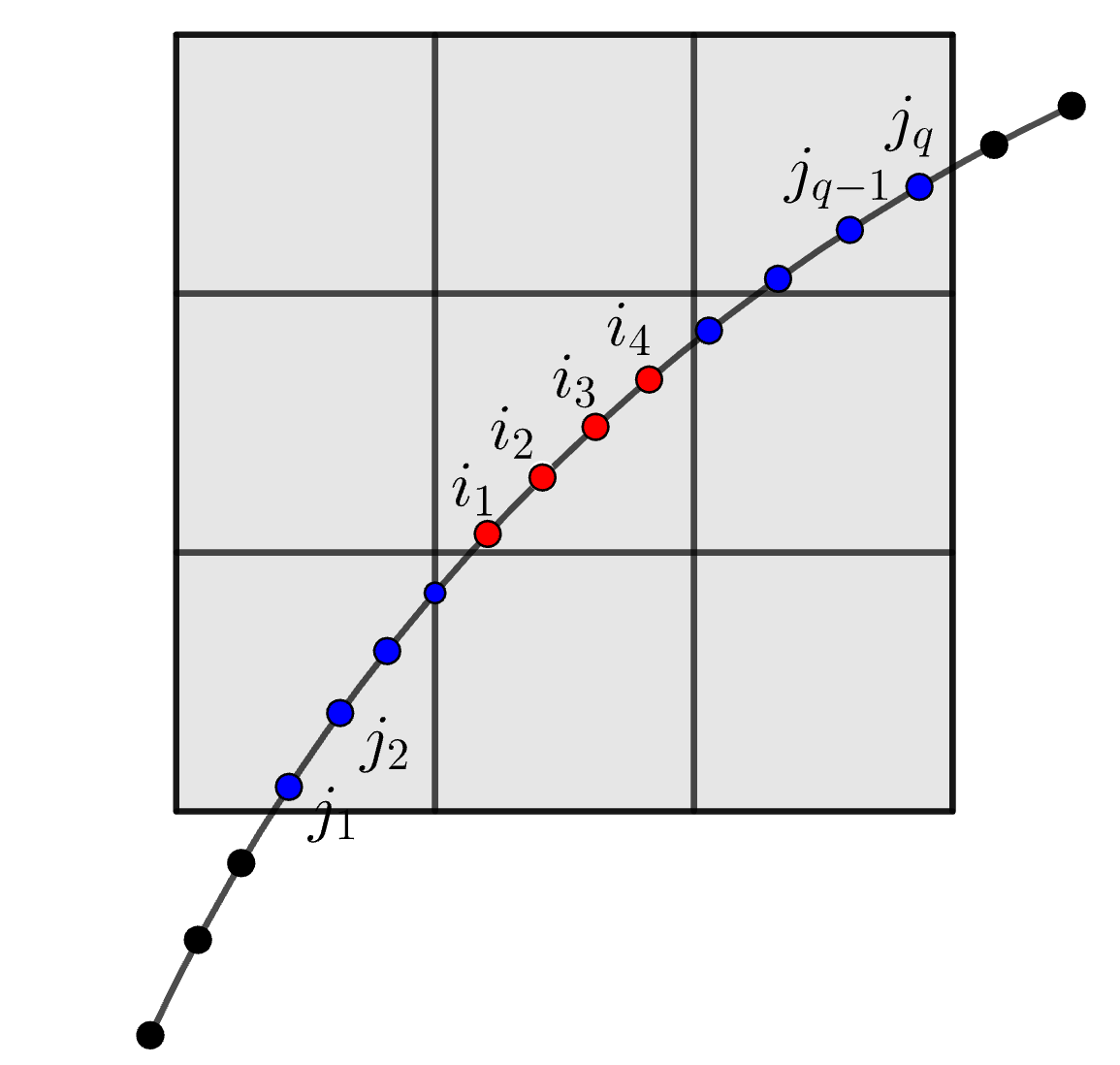}
\caption{Unknowns in the local system for preconditioning.}
\label{preconditioner}
\end{figure}

Within this FMM framework, source and target particles are the same. We now consider any given leaf box $B_{\ell k}$ in the source particle tree corresponding to the FMM for the free space component in the $\ell$-th layer. Let $B_{\ell k}$ contain particles associated with DOFs $\phi_{\ell i_1}, \phi_{\ell i_2}, \ldots, \phi_{\ell i_p}$, and let its neighbor boxes (including itself) contain DOFs $\phi_{\ell j_1}, \phi_{\ell j_2}, \ldots, \phi_{\ell j_q}$, see Fig. \ref{preconditioner}. Apparently, there holds
\[
\{i_1, i_2, \ldots, i_p\} \subset \{j_1, j_2, \ldots, j_q\}.
\]
We construct a local linear system
 \begin{equation*}
    {\mathbb{P}}_{B_{\ell k}}\begin{bmatrix}
        \tilde\phi_{\ell j_1}\\
        \tilde\phi_{\ell j_2}\\
        \vdots\\
        \tilde\phi_{\ell j_q}
    \end{bmatrix}=\begin{bmatrix}
        \eta_{\ell}g(\bm c_{\ell j_1})\\
        \eta_{\ell}g(\bm c_{\ell j_2})\\
        \vdots\\
        \eta_{\ell}g(\bm c_{\ell j_q})
    \end{bmatrix},
\end{equation*}
where 
\begin{equation} \label{localmatix}
   {\mathbb{P}}_{B_{\ell k}} = \begin{bmatrix}
        I_{\ell j_1}^s(\bm{c}_{\ell j_1}) & |\Gamma_{\ell j_2}|G_{\ell}^f(\bs c_{\ell j_2}, \bs c_{\ell j_1}) & \cdots & |\Gamma_{\ell j_q}|G_{\ell}^f(\bs c_{\ell j_q}, \bs c_{\ell j_1})\\
        |\Gamma_{\ell j_1}|G_{\ell}^f(\bs c_{\ell j_1}, \bs c_{\ell j_2}) & I_{\ell j_2}^s(\bm{c}_{\ell j_2}) & \cdots & |\Gamma_{\ell j_q}|G_{\ell}^f(\bs c_{\ell j_q}, \bs c_{\ell j_2})\\
        \vdots & \vdots & \ddots & \vdots\\
        | \Gamma_{\ell j_1}|G_{\ell}^f(\bs c_{\ell j_1}, \bs c_{\ell j_q}) & | \Gamma_{\ell j_2}|G_{\ell}^f(\bs c_{\ell j_2}, \bs c_{\ell j_q}) & \cdots   & I_{\ell j_q}^s(\bm{c}_{\ell j_q})
    \end{bmatrix}.
\end{equation}
Then, the components \(\tilde{\phi}_{\ell i_1}, \tilde{\phi}_{\ell i_2}, \ldots, \tilde{\phi}_{\ell i_p}\) corresponding to the particles within $B_{\ell k}$ are extracted from the solution as an approximation to the original unknowns \(\phi_{\ell i_1}, \phi_{\ell i_2}, \ldots, \phi_{\ell i_p}\). The inverse operator associated with the superposition of  these local solution procedures is denoted by \(\mathbb{P}_{\ell\ell}^{-1}\), satisfying:
\begin{equation*}
    \begin{bmatrix}
        \tilde\phi_{\ell 1}\\
        \tilde\phi_{\ell 2}\\
        \vdots\\
        \tilde\phi_{\ell N_{\ell}}
    \end{bmatrix}=\mathbb{P}_{\ell\ell}^{-1}\begin{bmatrix}
        \eta_{\ell}g(\bm c_{\ell 1})\\
        \eta_{\ell}g(\bm c_{\ell 2})\\
        \vdots\\
        \eta_{\ell}g(\bm c_{\ell N_{\ell}})
    \end{bmatrix}.
\end{equation*}
Suppose the FMM tree structure for the \(\ell\)-th layer has \(m_{\ell}\) leaf boxes, then
\[\mathbb{P}_{\ell \ell}^{-1} = \sum_{k=1}^{m_{\ell}}\mathbb{S}_{B_{\ell k}}\mathbb{Q}_{B_{\ell k}}\mathbb{P}_{B_{\ell k}}^{-1}\mathbb{R}_{B_{\ell k}},\]
where $\mathbb{R}_{B_{\ell k}}$, $\mathbb{Q}_{B_{\ell k}}$ and $\mathbb{S}_{B_{\ell k}}$ are gather and scatter matrices corresponding to box $B_{\ell k}$ such that 
$$\begin{bmatrix}
        \phi_{\ell j_1}\\
        \phi_{\ell j_2}\\
        \vdots\\
        \phi_{\ell j_q}
    \end{bmatrix}
    =\mathbb R_{B_{\ell k}}{\bf\Phi}_{\ell},\quad \begin{bmatrix}
        \phi_{\ell i_1}\\
        \phi_{\ell i_2}\\
        \vdots\\
        \phi_{\ell i_p}
    \end{bmatrix}=\mathbb Q_{B_{\ell k}}\begin{bmatrix}
        \phi_{\ell j_1}\\
        \phi_{\ell j_2}\\
        \vdots\\
        \phi_{\ell j_q}
    \end{bmatrix},
    \quad {\bf\Phi}_{\ell}=\sum_{k=1}^{m_{\ell}}\mathbb S_{B_{\ell k}}\begin{bmatrix}
        \phi_{\ell i_1}\\
        \phi_{\ell i_2}\\
        \vdots\\
        \phi_{\ell i_p}
    \end{bmatrix}.$$
Consequently, 
\begin{equation*}
{\mathbb P}^{-1}=   \begin{bmatrix}
     {\mathbb P}_{00}^{-1} & 0 & \cdots & 0\\
     0 & {\mathbb P}_{11}^{-1} & \cdots & 0\\
     \vdots &\vdots & \ddots & \vdots &\\
     0 & 0 & \cdots & {\mathbb P}_{LL}^{-1}
    \end{bmatrix}
\end{equation*}
completes the construction of a high-performance preconditioner for the original matrix \(\mathbb{K}\). The effectiveness of the proposed preconditioner will be validated by numerical examples presented in the next section.

\section{Numerical Examples}\label{sect5}

In this section, we present some numerical examples to validate the accuracy and efficiency of the proposed fast boundary integral method to solve multiple scattering problems in layered media. All examples are solved using the lowest order Nystr\"{o}m discretization and FMM accelerated GMRES algorithm with a relative error tolerance of \(1.0e-8\).

\begin{figure}[htbp] 
    \centering 
    \subfigure[layered media diagram]{\includegraphics[width=0.32\textwidth]{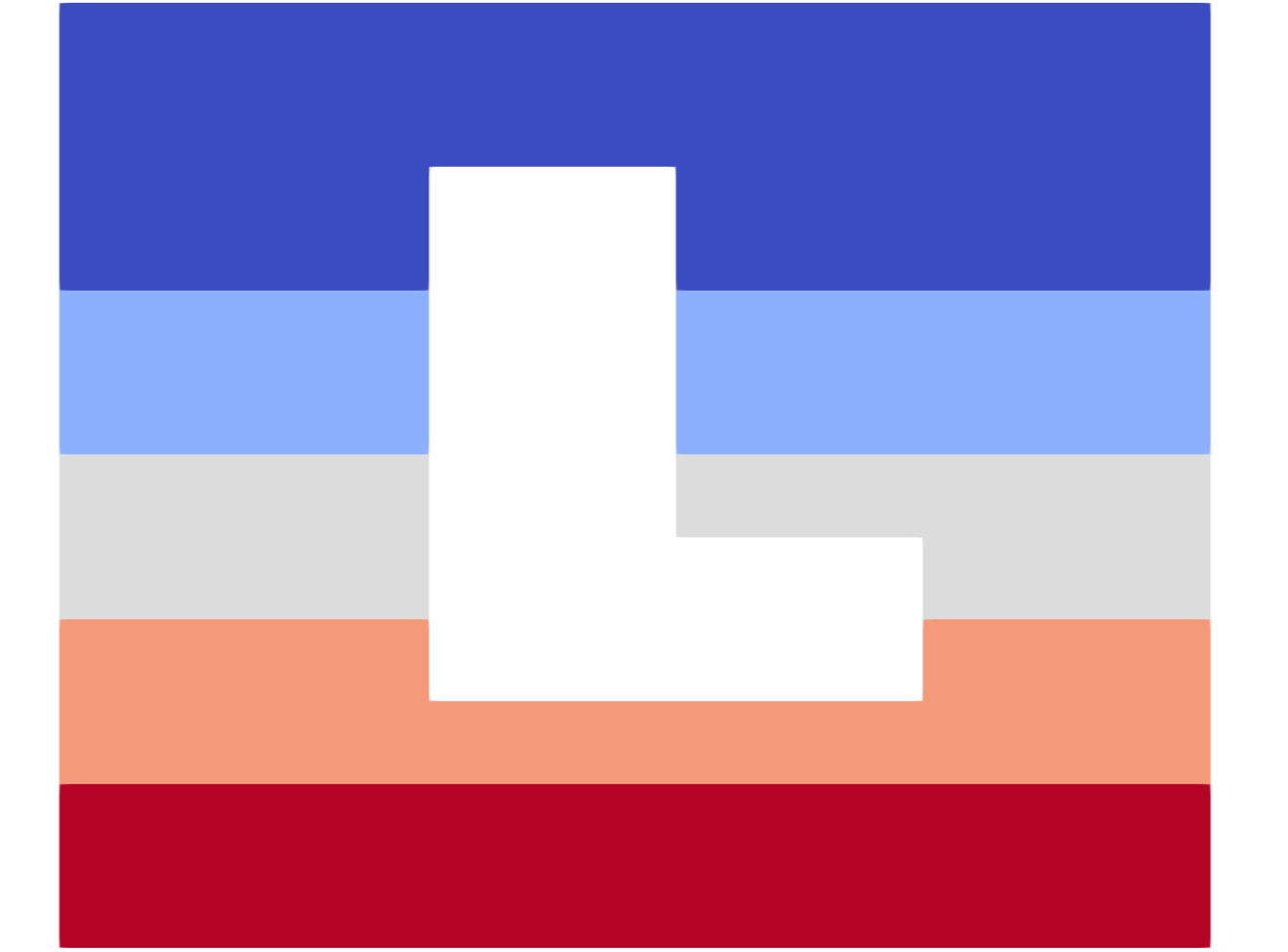}}
    \subfigure[$\mathfrak{Re}(u^s)$]{\includegraphics[width=0.32\textwidth]{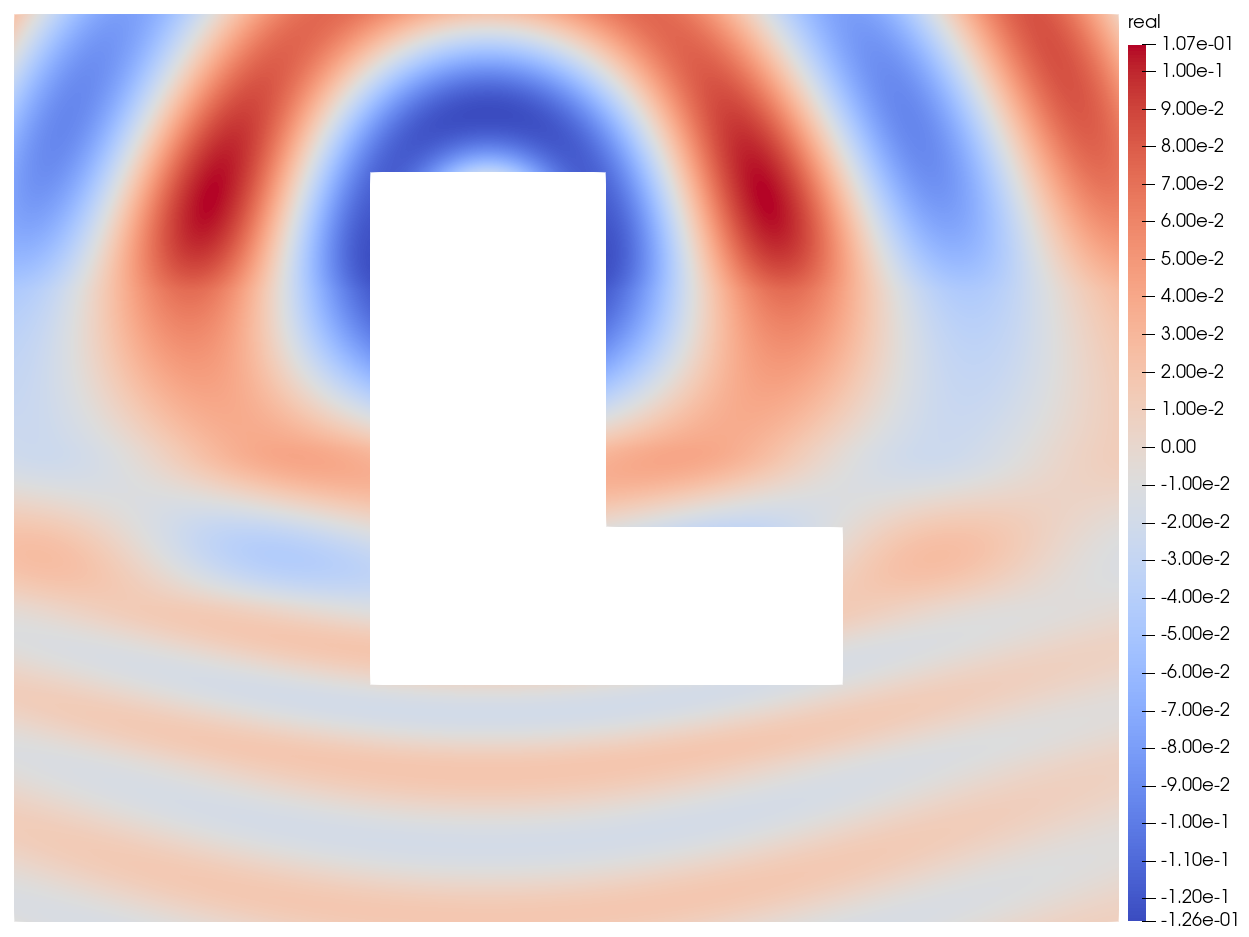}}
    \subfigure[$\mathfrak{Im}(u^s)$]{\includegraphics[width=0.32\textwidth]{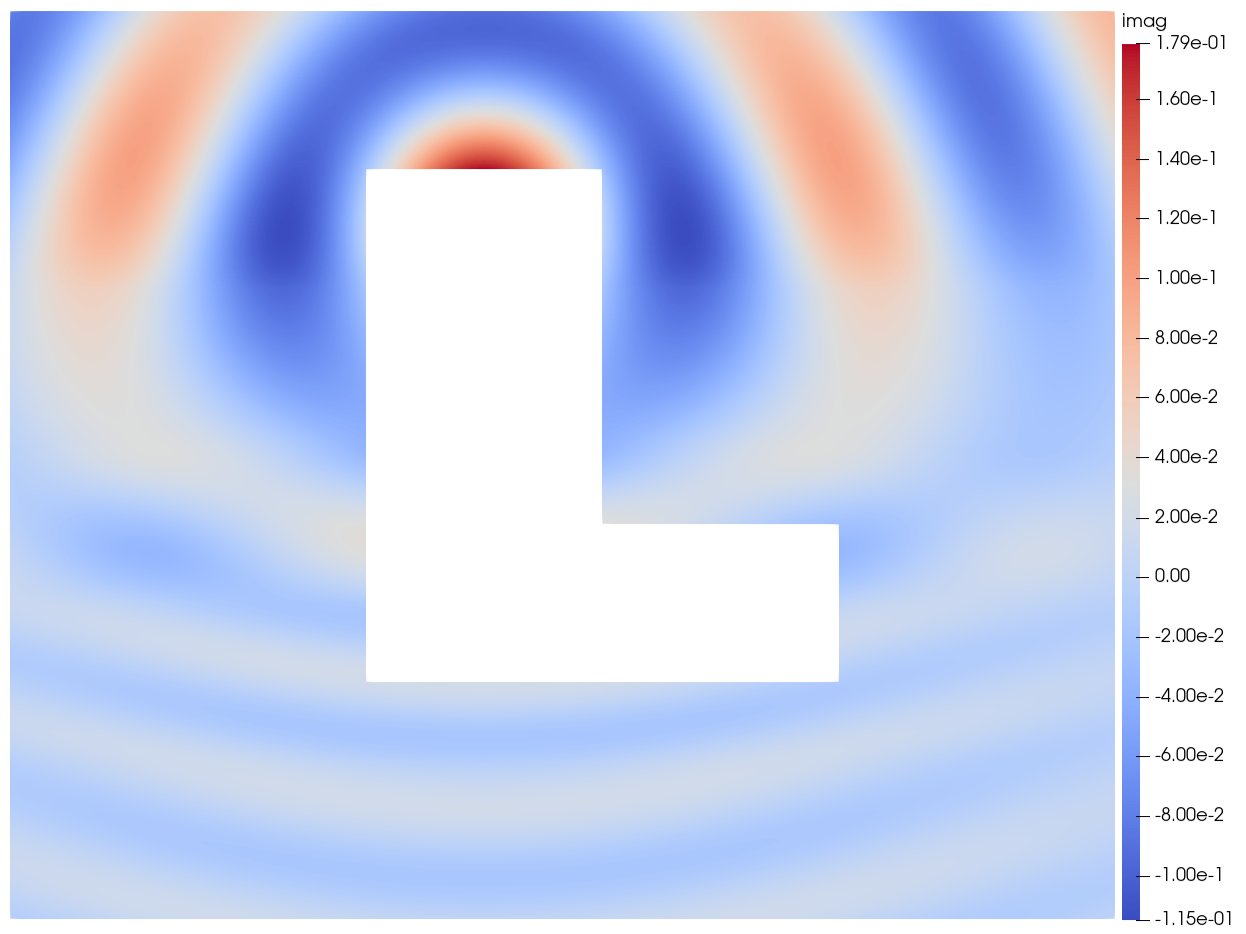}}
    \caption{Problem diagram and numerical solution} 
    \label{Lshape} 
\end{figure}

\vspace{10pt}
\noindent \textbf{Example 1.} We first test the accuracy of the proposed fast algorithm to solve the scattering problem \eqref{extdirichleteq1} with an exact solution given by
\begin{equation}\label{numericalexamples}
    u^{\rm{s}}(\bm{r})=G_{\ell 0}(\bm{r},\bm{r}'),\quad \forall\bm{r}\in \Upsilon_{\ell}\setminus \Omega ,\quad \ell=0, 1, 2, 3, 4,
\end{equation}
where $\Omega$ is a L-shaped domain with vertices (0.75, 0.75), (-0.75, 0.75), (-0.75, -2.5), (2.25, -2.5), (2.25, -1.5), and (0.75, -1.5), as shown in Fig. \ref{Lshape}(a), the interfaces of the layered medium are located at $y=\{0,-1,-2,-3\}$. Apparently, the scatterer $\Omega$ is embedded in 4 layers. The reflection indices and the wave numbers of the layered media are specified as $\eta = \{1, 2, 3, 4, 5\}$ and $k = \{3.2, 2.5, 5.1, 8.6, 6.9\}$ . The exact solution \(G_{\ell 0}\) is the Green's function of the layered media with a point source at $\bs r'=(0, 0.375)$ located inside the L-shaped domain in the $0$-th layer.

\begin{figure}[H]
    \centering
    \captionsetup{skip=5pt} 
    \begin{minipage}[b]{0.48\textwidth}
        \centering
        \includegraphics[width=\textwidth]{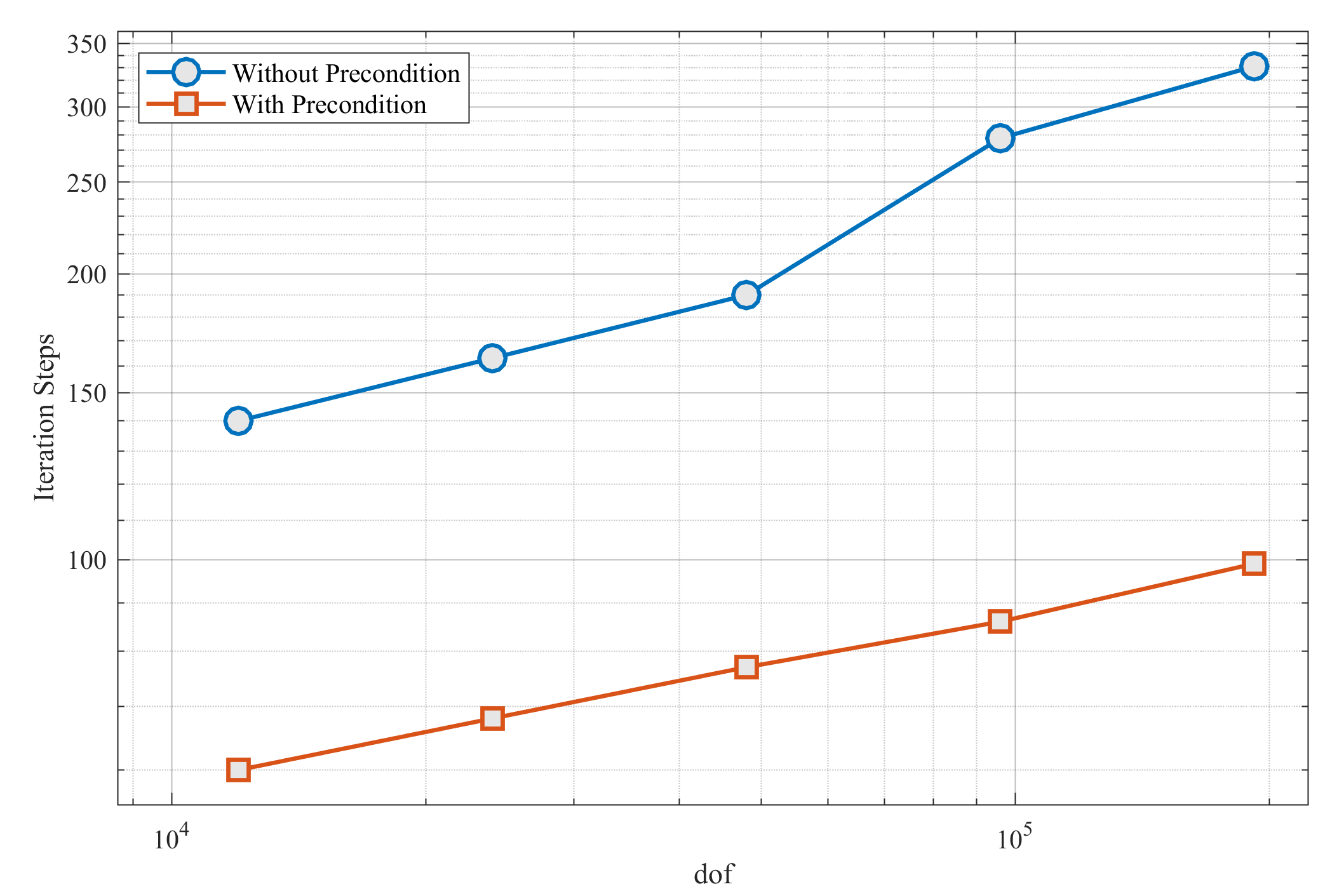}
        \caption{Iteration steps vs dof}
        \label{Iteration_steps}
    \end{minipage}
    \hfill
    \begin{minipage}[b]{0.48\textwidth}
        \centering
        \includegraphics[width=\textwidth]{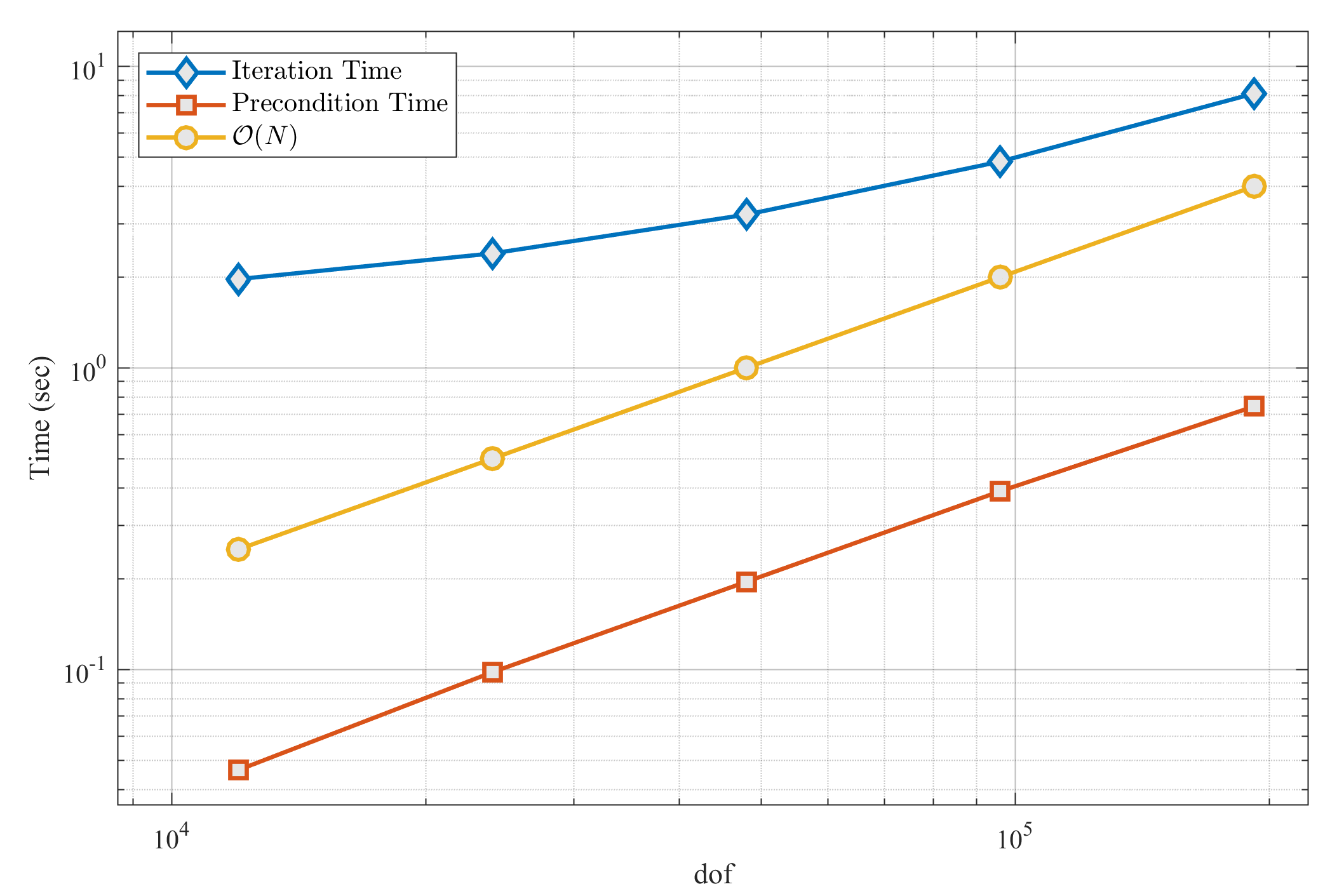}
        \caption{Time per iteration step vs dof}
        \label{Iteration_time}
    \end{minipage}
\end{figure}

The performance of the proposed fast algorithm is summarized in Table \ref{example1NumericalResults}. It shows that first-order convergence rates are achieved in both the \(L_{\infty}\) and \(L_2\) norms. By using the preconditioning, the number of iterations decreases significantly, especially when the degree of freedom is large. The numbers of iterations in the GMRES with/without preconditioning are compared in Table \ref{example1NumericalResults} and Fig. \ref{Iteration_steps}. Besides, the CPU time of the FMM accelerated iteration and the preconditioning in each step are plotted in Fig. \ref{Iteration_time}. Both show \(\mathcal{O}(N)\) computational complexity. 
\begin{table}[H]
    \centering
    \caption{The performance of the fast algorithm}
    \label{example1NumericalResults}
    \scalebox{0.8}{ 
    \begin{tabular}{ccccccc}
        \toprule 
        Elements & Iteration Steps & \makecell{Iteration Steps \\ with preconditioning} & \makecell{Time(sec)} & \makecell{Time(sec)  with \\ preconditioning} & $L_{\infty}$ Error & $L_{2}$ Error \\
        \midrule 
        12000 & 140 &60 & 739.71 &573.37 & 2.282e-02 & 3.064e-04 \\
        24000 & 163 &68 & 918.03 &650.14 & 1.151e-02 & 1.544e-05 \\
        48000 & 190 &77 & 1142.87 &783.46 & 5.782e-03 & 7.748e-05 \\
        96000 & 278 &86 & 2110.27 &1014.66 & 2.897e-03 & 3.884e-05\\
        192000 & 331 &99 & 3535.53 &1501.94 & 1.450e-03 & 1.948e-05\\
        \bottomrule 
    \end{tabular}
    }
\end{table}

\begin{figure}[htbp] 
    \centering 
    \subfigure[layered media diagram]{\includegraphics[width=0.32\textwidth]{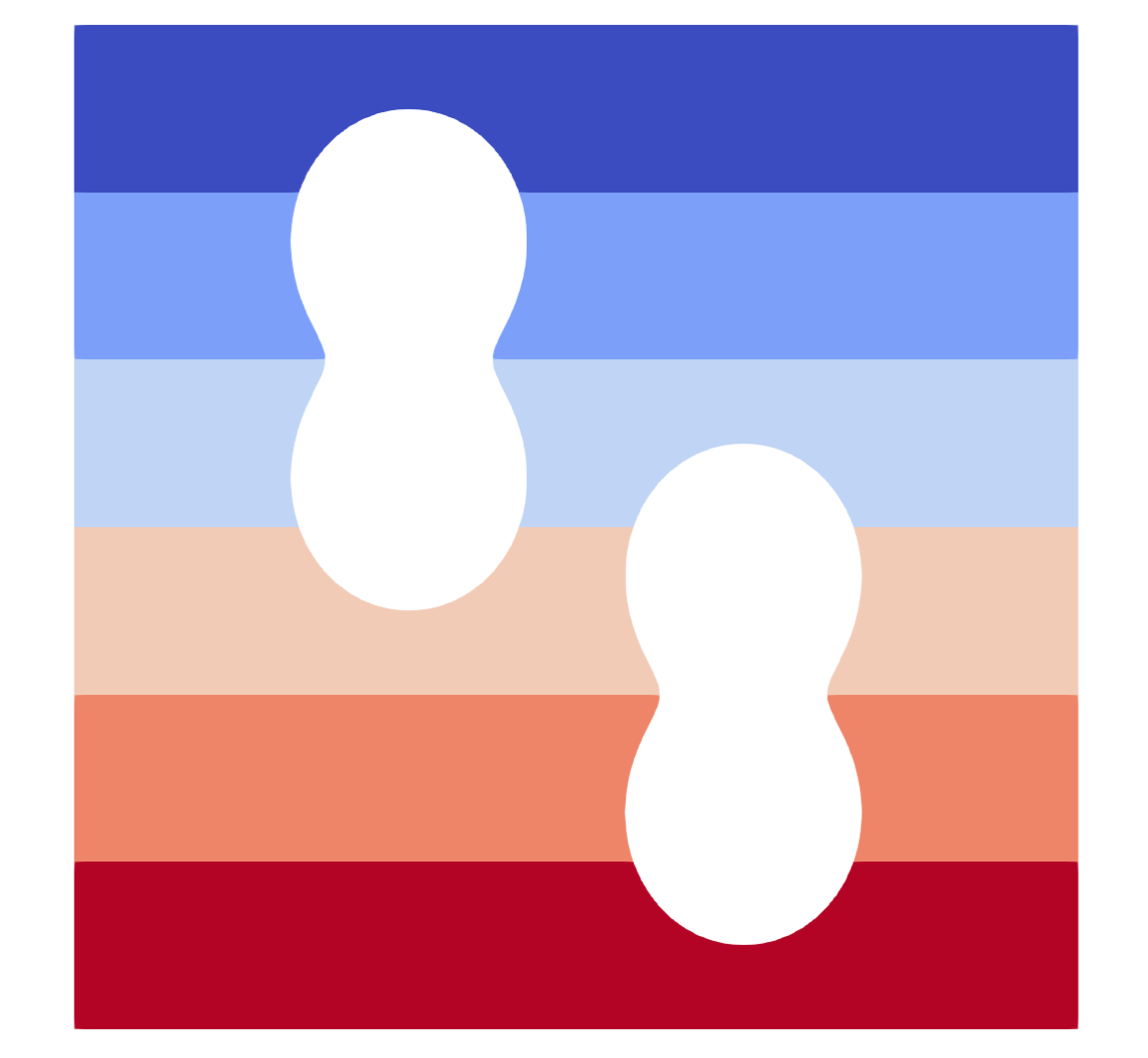}}
    \subfigure[$\mathfrak{Re}(u^s)$]{\includegraphics[width=0.32\textwidth]{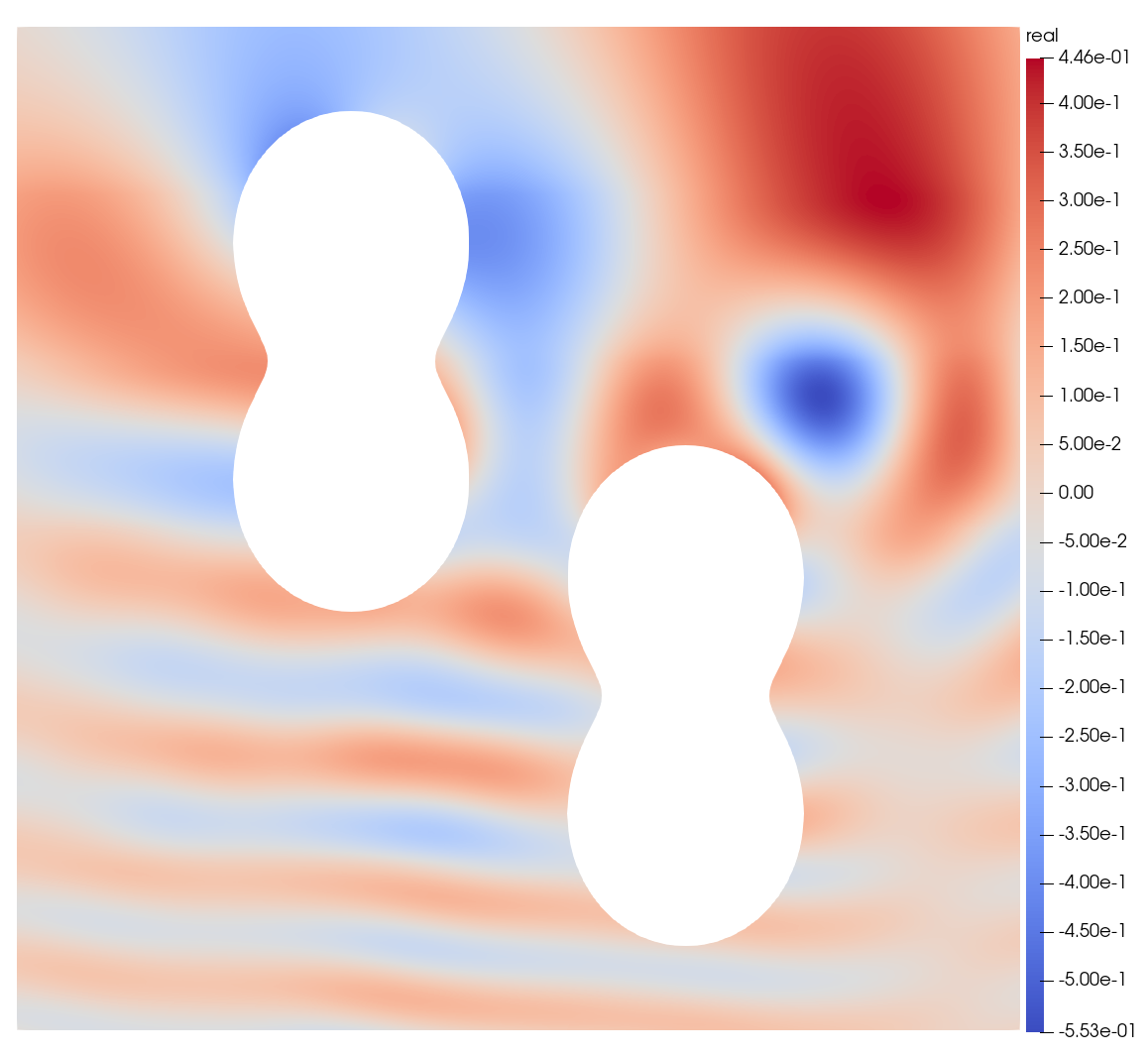}}
    \subfigure[$\mathfrak{Im}(u^s)$]{\includegraphics[width=0.32\textwidth]{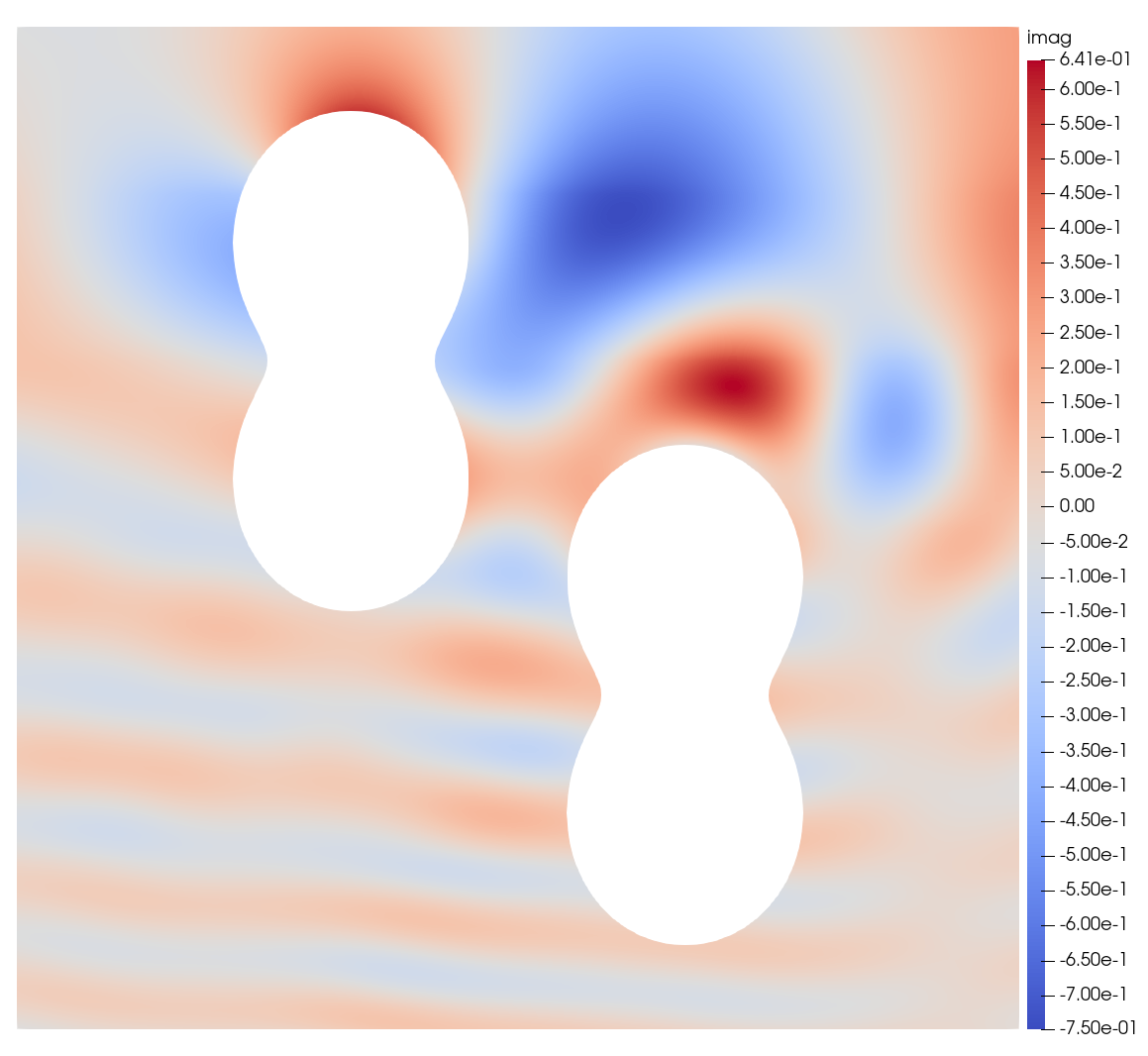}}
    \caption{Problem diagram and numerical solution} 
    \label{pernut} 
\end{figure}

\begin{figure}[H]
    \centering
    \captionsetup{skip=5pt} 
    \begin{minipage}[b]{0.48\textwidth}
        \centering
        \includegraphics[width=\textwidth]{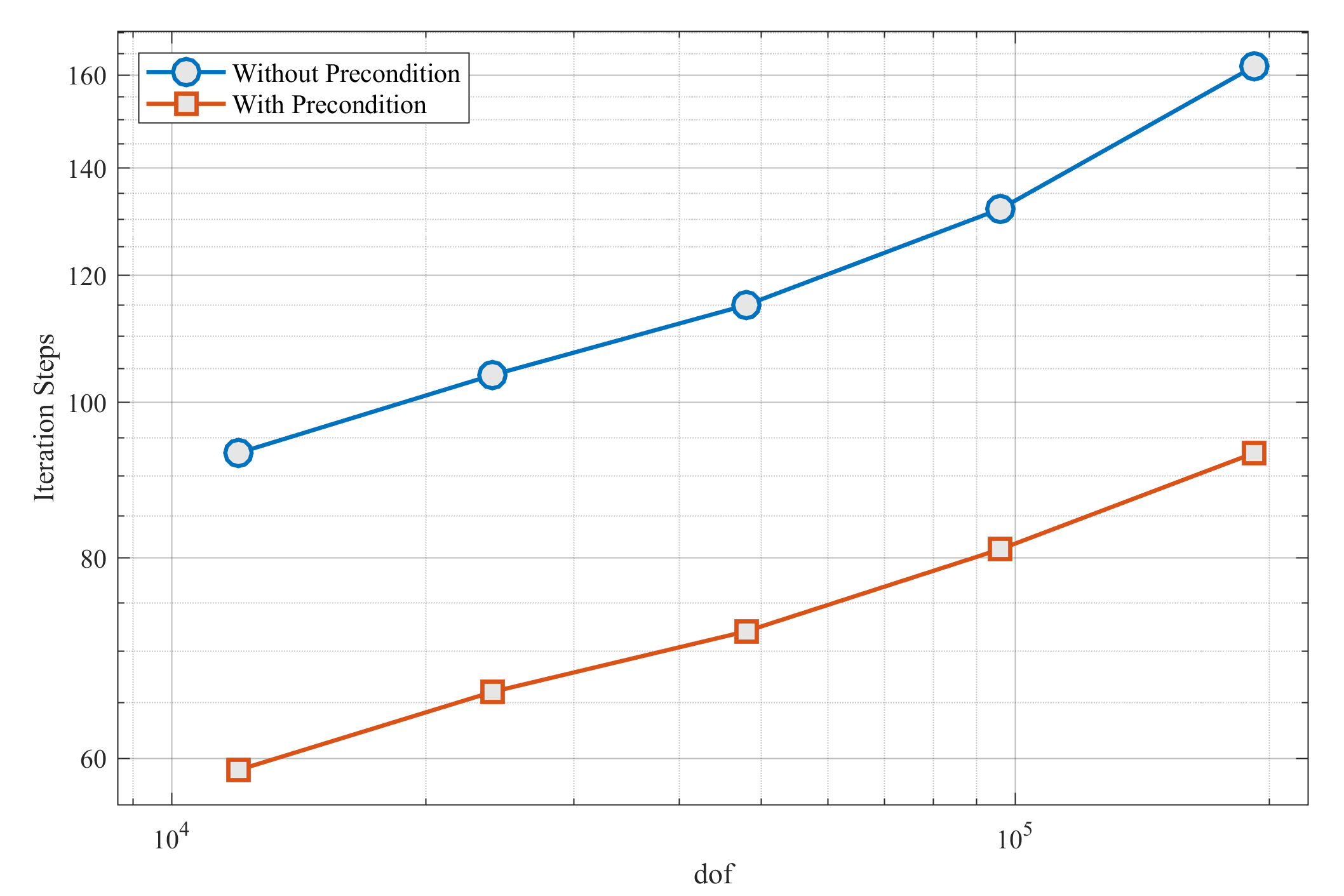}
        \caption{Iteration steps vs dof}
        \label{pernutIterationsteps}
    \end{minipage}
    \hfill
    \begin{minipage}[b]{0.48\textwidth}
        \centering
        \includegraphics[width=\textwidth]{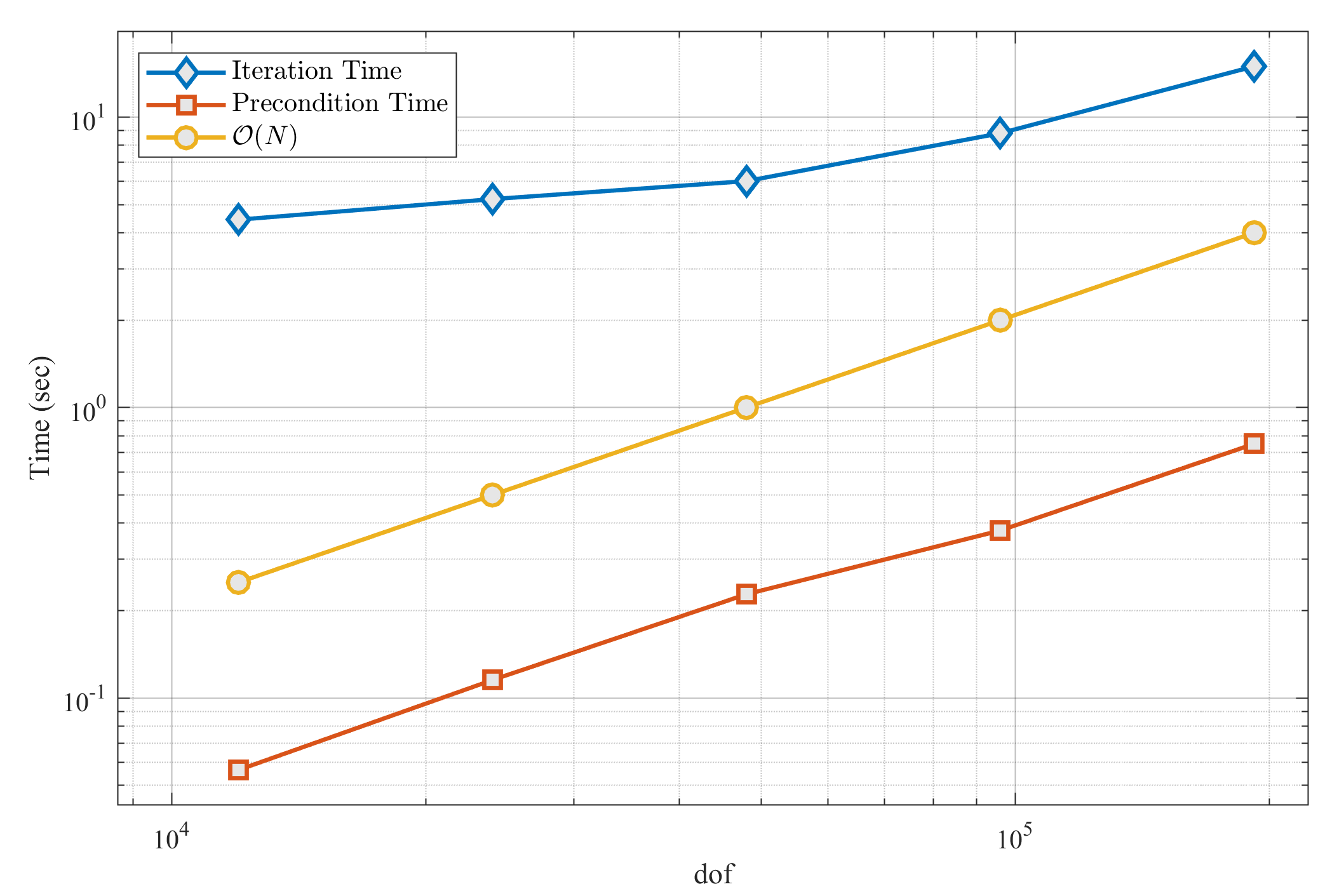}
        \caption{Time per iteration step vs dof}
        \label{pernutIterationtime}
    \end{minipage}
\end{figure}

\noindent \textbf{Example 2.} Consider a layered medium of six layers with interfaces at $y=\{0,-1,-2,-3,-4\}$, where the reflection indices and wave numbers are specified as $\eta = \{1.1, 2.3, 3.4, 4.6, 5.0, 6.6\}$ and $k = \{1.2,2.3,4.5,6.1,7.7,10.0\}$. We shall use  parametrics
\begin{equation}\label{scattererpara}
r_i=a\sin k(\theta_i-\theta_0)+b,\quad \theta\in [0, 2\pi],
\end{equation}
to generate scatterers, where $(r_i,\theta_i)$ is the polar coordinates of $\bs r$ with respect to given center $\bs c_i$.
Here, two scatterers (see Fig. \ref{pernut}(a)) generated by (\ref{scattererpara}) with $k=2, a=0.5, b=1, \theta_0 = \frac{\pi}{4}$ and $\bs{c}_1=(0,-1)$ and $\bs{c}_2=(2,-3)$ are tested.  Given incident wave \(u^{\rm{inc}} = e^{\ri(-\frac{3\sqrt{2}}{5}x-\frac{3\sqrt{2}}{5}y)}\), the scattering field  \(u^s\) are plotted in Figs. \ref{pernut}(b) (real part) and \ref{pernut}(c) (imaginary part). 
In Table \ref{pernuIterationsdata}, the number of iterations and the total CPU time are presented. Further, the number of iterations and CPU time of the FMM accelerated iteration and the preconditioning in each iteration are plotted in Figs. \ref{pernutIterationsteps} and \ref{pernutIterationtime}. The results show that the preconditioning is very effective and  the computational complexity of each iteration step is \(\mathcal{O}(N)\). 
\begin{table}[htbp]
    \centering
    \caption{Iteration Steps and Iteration time}
    \label{pernuIterationsdata}
    \begin{tabular}{ccccc}
        \toprule 
        Elements & Iteration Steps & \makecell{Iteration Steps \\ with preconditioning} & \makecell{Time(sec)} & \makecell{Time(sec)  with \\ preconditioning} \\
        \midrule 
        12000 & 93 &59 &1628.33 & 1362.33  \\
        24000 & 104 &66 &1822.82 &1520.70  \\
        48000 & 115 &72 &2429.23 &1685.44  \\
        96000 & 132 &81 &3042.68 &2058.86\\
        192000 & 162 &93 &3807.98 &2827.82\\
        \bottomrule 
    \end{tabular}
\end{table}

\begin{figure}[htbp] 
    \centering 
    \subfigure[layered media diagram]{\includegraphics[width=0.32\textwidth]{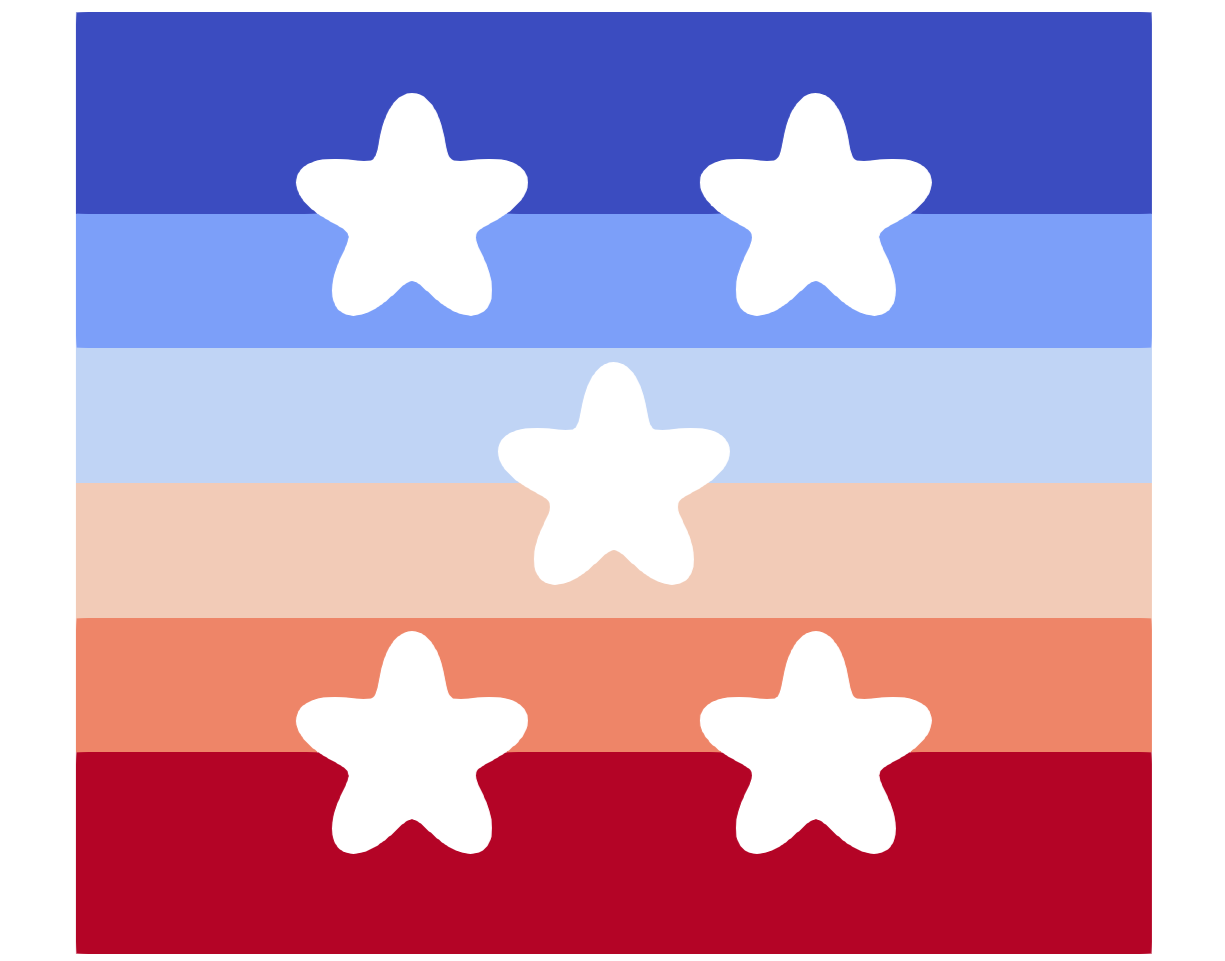}}
    \subfigure[$\mathfrak{Re}(u^s)$]{\includegraphics[width=0.32\textwidth]{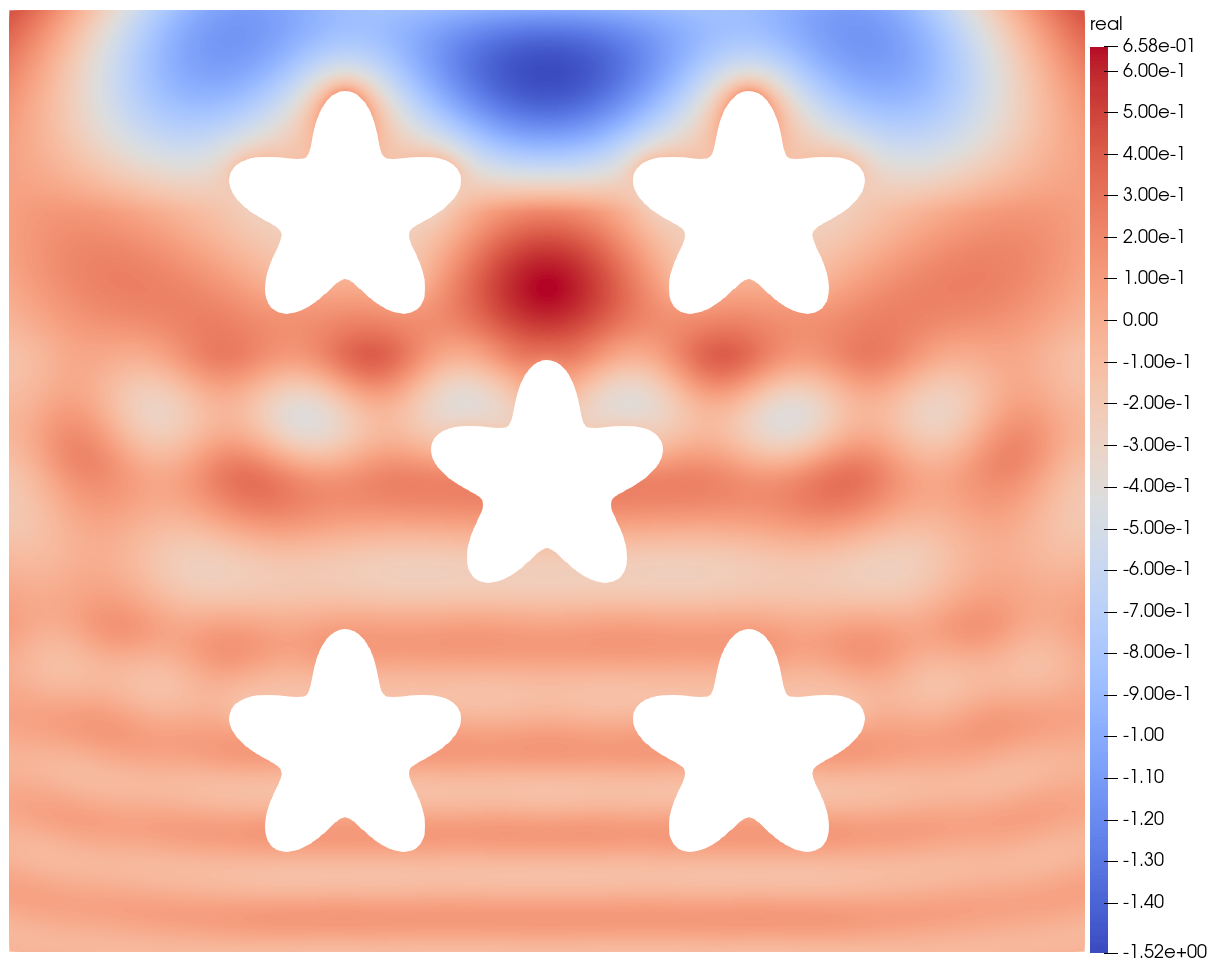}}
    \subfigure[$\mathfrak{Im}(u^s)$]{\includegraphics[width=0.32\textwidth]{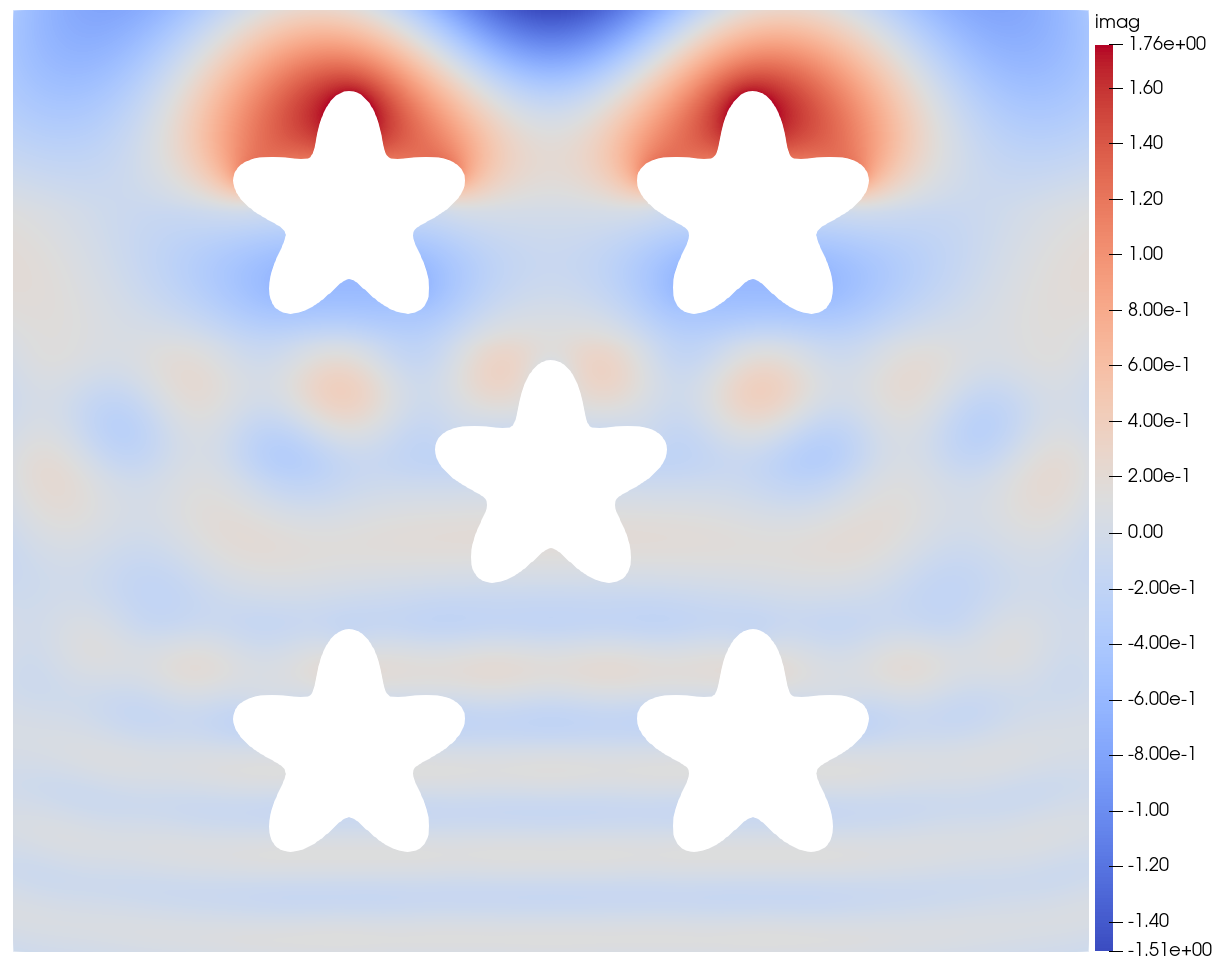}}
    \caption{Problem diagram and numerical solution} 
    \label{flower} 
\end{figure}

\begin{figure}[H]
    \centering
    \captionsetup{skip=5pt} 
    \begin{minipage}[b]{0.48\textwidth}
        \centering
        \includegraphics[width=\textwidth]{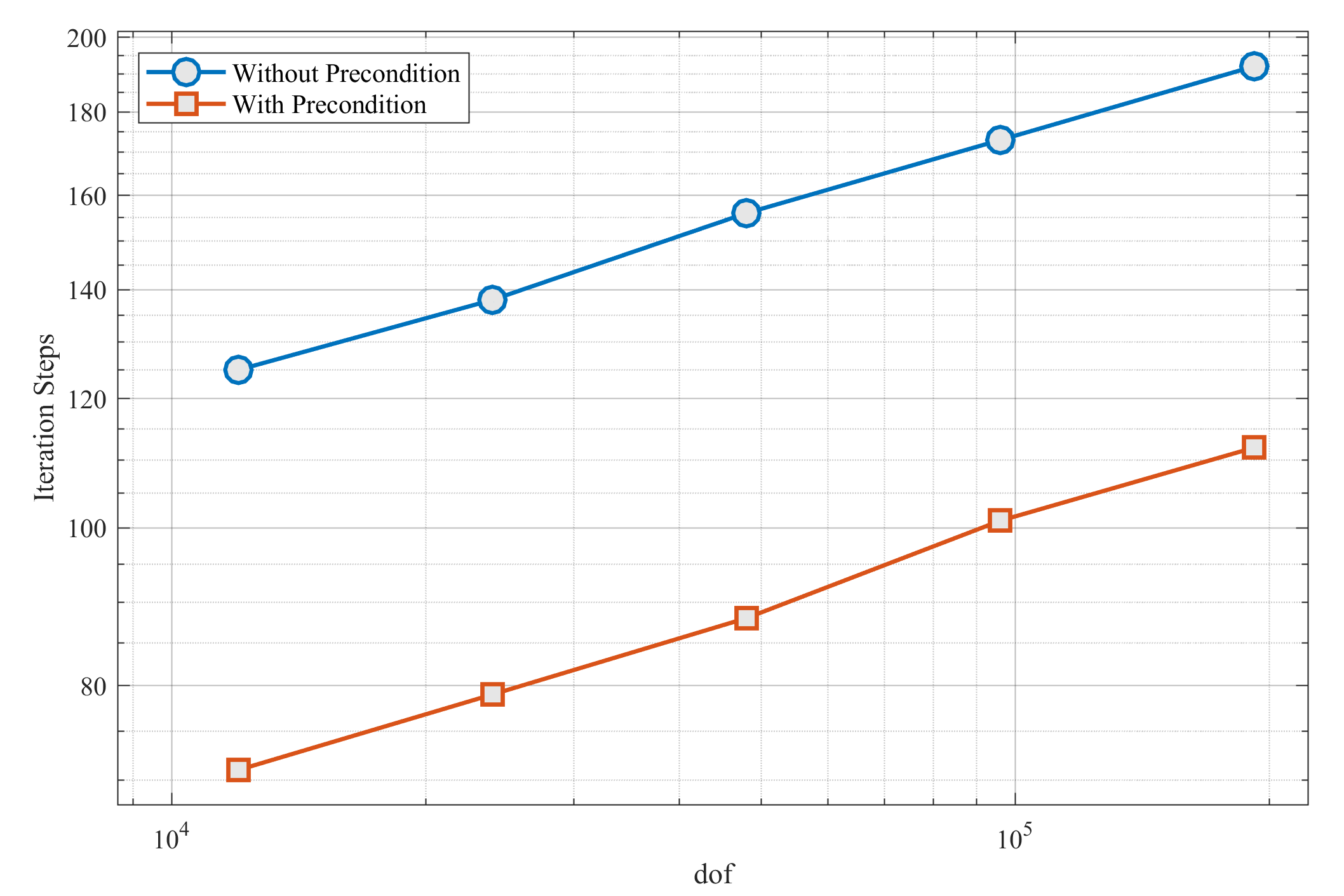}
        \caption{Iteration steps vs dof}
        \label{flowerIterationsteps}
    \end{minipage}
    \hfill
    \begin{minipage}[b]{0.48\textwidth}
        \centering
        \includegraphics[width=\textwidth]{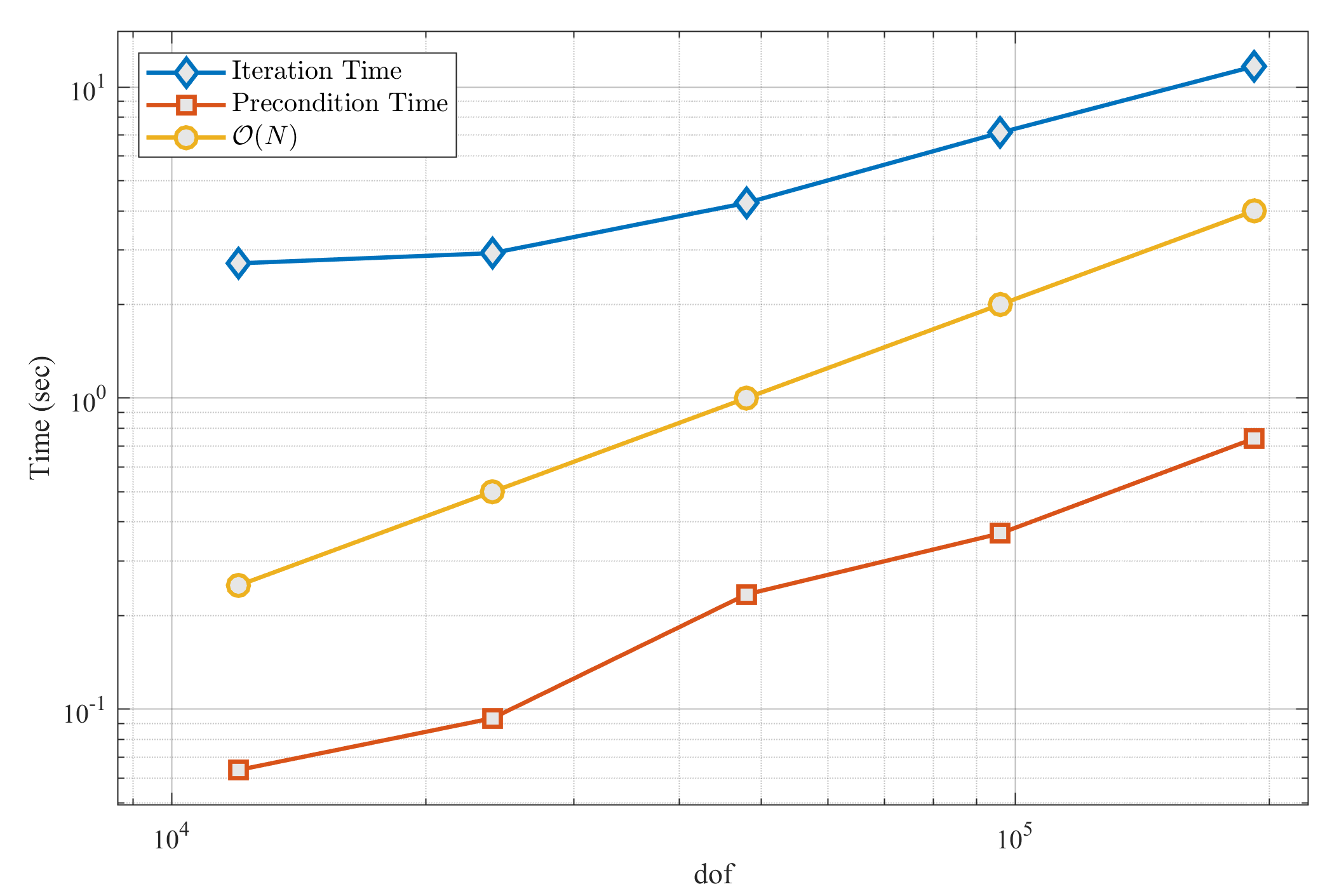}
        \caption{Time per iteration step vs dof}
        \label{flowerIterationtime}
    \end{minipage}
\end{figure}

\noindent \textbf{Example 3.} Consider the layered medium of six layers with interfaces at $y=\{0,-1,-2,-3,-4\}$, where the reflection indices and the wave numbers are specified as $\eta = \{1, 2, 3, 4, 5,6\}$ and $k = \{2, 3, 6, 5, 8 ,10\}$. 
Five scatterers (see Fig. \ref{flower}(a)) determined by (\ref{scattererpara}) with $k=5, a=0.2, b=0.7, \theta_0 = 0$ and $\bs{c}_1=(-1.5,0)$, $\bs{c}_2=(1.5,0)$, $\bs{c}_3=(0,-2)$, $\bs{c}_4=(-1.5,-4)$ and $\bs{c}_5=(1.5,4)$ are tested. Given incident wave \(u^{\rm{inc}} = e^{-\ri2y}\), the scattering field \(u^s\) are plotted in Figs. \ref{flower}(b) (real part) and \ref{flower}(c) (imaginary part).
The performance of the fast algorithm is presented in Table \ref{flowerIterationsdata}, Fig. \ref{flowerIterationsteps} and Fig. \ref{flowerIterationtime}. Apparently, the same conclusion on the efficiency of the algorithm can be obtained as we have validated in Example 2. 

\begin{table}[htbp]
    \centering
    \caption{Iteration Steps and Iteration time}
    \label{flowerIterationsdata}
    \scalebox{1.0}{ 
    \begin{tabular}{ccccc}
        \toprule 
        Elements & Iteration Steps & \makecell{Iteration Steps \\ with preconditioning} & \makecell{Time(sec)} & \makecell{Time(sec)  with \\ preconditioning} \\
        \midrule 
        12000 & 125 &71 & 1638.88 &1553.62  \\
        24000 & 138 &79 & 1742.52 &1639.49  \\
        48000 & 156 &88 & 2031.50 &1831.76 \\
        96000 & 173 &101 & 2708.44 &2649.50 \\
        192000 & 192 &112 & 3782.74 &3056.38 \\
        \bottomrule
    \end{tabular}
    }
\end{table}

\vskip 20pt
\noindent{\large\bf Conclusion and future work}
\vskip 10pt

In this paper, we have developed a FMM accelerated boundary integral method for solving two-dimensional Helmholtz equation for acoustic wave scattering problems in multi-layered media. By using the layered Green’s functions to formulate the boundary integral equations, the interface transmission conditions are naturally satisfied without redundant unknowns and integral equations on the interfaces. The lowest order Nystr\"{o}m method is used to discretize the integral equation. An improved layered media FMM and an overlapped domain decomposition preconditioner are proposed to accelerate the matrix free GMRES iterative solver for the resulted dense linear system. A set of numerical experiments demonstrated the $\mathcal O(N)$ complexity and robustness of the method, confirming its capability for efficient scattering simulations in complex multilayer configurations.

For the future work, the extension of the present framework to three-dimensional layered media will be considered, where the complexity of Sommerfeld integrals and interface coupling poses additional challenges. Besides, further advances in preconditioning strategies and parallel implementations will be considered to enable large-scale simulations in high-frequency regimes. 


\vskip 10pt
\noindent{\large\bf Acknowledgment}
\vskip 6pt
The research of the first and second author is  supported by NSFC (Grant No. 12022104 and 12371394). The research of the third author is support by the Major Program of Xiangjiang Laboratory(Grant No. 22XJ01013) and the Foundation Sciences of Hunan Province (Grant No. 2024JC0006). The research of W.  Cai is supported by Clements Chair of Applied Mathematics at Southern Methodist University.


\end{document}